\documentclass[sn-mathphys,Numbered]{sn-jnl}

\usepackage{graphicx}%
\usepackage{multirow}%
\usepackage{amsmath,amssymb,amsfonts}%
\usepackage{amsthm}%
\usepackage{bm}
\usepackage{mathtools}
\usepackage{mathrsfs}%
\usepackage[title]{appendix}%
\usepackage{xcolor}%
\usepackage{textcomp}%
\usepackage{manyfoot}%
\usepackage{booktabs}%
\usepackage{algorithm}%
\usepackage{algorithmicx}%
\usepackage{algpseudocode}%
\usepackage{listings}%
\usepackage{caption}
\usepackage{subcaption}
\newcommand{\imag}{\texttt{i}}

\begin{document}

\title[Article Title]{Dynamics Around the Earth-Moon Triangular Points in the Hill Restricted 4-Body Problem}


\author*[1]{\fnm{Luke T.} \sur{Peterson}}\email{Luke.Peterson@colorado.edu}

\author[1]{\fnm{Gavin} \sur{Brown}}\email{gavin.m.brown@colorado.edu}

\author[2,3]{\fnm{\`Angel} \sur{Jorba}}\email{angel@maia.ub.es}

\author[1]{\fnm{Daniel} \sur{Scheeres}}\email{scheeres@colorado.edu}

\affil*[1]{\orgdiv{Ann and H.J. Smead Aerospace Engineering Sciences}, \orgname{University of Colorado Boulder}, \state{CO}, \country{USA}}

\affil[2]{\orgdiv{Department de Matem\`atiques i Inform\`atica}, \orgname{Universitat de Barcelona (UB)}, \city{Barcelona}, \country{Spain}}

\affil[3]{\orgdiv{Centre de Recerca Matem\`atica}, \orgname{Edifici C}, \orgaddress{\street{Campus Bellaterra}, \postcode{08193} \city{Bellaterra}, \country{Spain}}}

\abstract{This paper investigates the motion of a small particle moving near the triangular points of the Earth-Moon system. The dynamics are modeled in the Hill restricted 4-body problem (HR4BP), which includes the effect of the Earth and Moon as in the circular restricted 3-body problem (CR3BP), as well as the direct and indirect effect of the Sun as a periodic time-dependent perturbation of the CR3BP. Due to the periodic perturbation, the triangular points of the CR3BP are no longer equilibrium solutions; rather, the triangular points are replaced by periodic orbits with the same period as the perturbation. Additionally, there is a 2:1 resonant periodic orbit that persists from the CR3BP into the HR4BP. In this work, we investigate the dynamics around these invariant objects by performing a center manifold reduction and computing families of 2-dimensional invariant tori and their linear normal behavior. We identify bifurcations and relationships between families. Mechanisms for transport between the Earth, $L_4$, and the Moon are discussed. Comparisons are made between the results presented here and in the bicircular problem (BCP).}


\maketitle

\newpage 

\tableofcontents

\newpage 

\section{Introduction}\label{sec:Introduction}

The study of dynamics in the Earth-Moon (EM) environment has become the focus of much research due to the announcement of the lunar gateway placed in a southern $L_2$ near rectilinear halo orbit (NRHO) \cite{smith2020artemis}, as well as nearby orbits used for Artemis program operations \cite{mccarthy2023four}. As the traffic of spacecraft increases for operations near the Moon, there is an increased need for relay communications between Earth and the vicinity of the Moon. TYCHO is a proposed mission to provide such communication relay services via satellites placed near the Earth-Moon $L_4$ libration point--one of the two so-called triangular points \cite{hornig2015tycho}. To model the dynamics in cislunar space, the circular restricted 3-body problem (CR3BP) is used as a first approximation \cite{szebehely_modifications_1967,gomez_dynamics_2001}. In the Earth-Moon CR3BP, the triangular points $L_{4,5}$ are elliptic equilibria--they are linearly stable. In fact, it has been shown that there is a region of effective stability around $L_4$ in the CR3BP \cite{giorgilli1989effective}. This stability, along with the fact that $L_4$ is equidistant to the Earth and Moon, presents an ideal region of cislunar space to put communications satellites.

The presence of the Trojan asteroids near the Sun-Jupiter triangular points has raised questions about the existence of similar objects in the Earth-Moon system \cite{gimeno2024effect}. After all, the triangular points are elliptic in both the Sun-Jupiter and Earth-Moon CR3BP, and Earth-Moon Trojans could pose a risk of collision with communications satellites. While it is known that there are no Earth-Moon Trojans, the existence of Kordylewski dust clouds around EM $L_4$ has recently been confirmed \cite{horvath2018kordylewski_1,horvath2018kordylewski_2}. The key difference between the Sun-Jupiter and Earth-Moon systems is that there are no other major perturbing forces in the Sun-Jupiter system, as the Sun and Jupiter are the two most massive celestial bodies in the solar system, while the Earth-Moon CR3BP neglects the significant perturbing force of the Sun's gravity. It has been shown that accounting for the gravitational effect of the Sun in the Earth-Moon system qualitatively changes the behavior of $L_4$--its stability changes from elliptic to partially hyperbolic \cite{scheeres1998restricted,simo1995bicircular}. Accounting for the Sun's gravity helps to describe the non-existence of Earth-Moon Trojan asteroids, and recent work has incorporated additional non-gravitational perturbations caused by the Sun to help describe the existence of the Kordylewski clouds from a dynamical astronomy perspective \cite{jorba2021stabilizing,gimeno2024effect}. However, what are the implications of this qualitative change in dynamical behavior around EM $L_4$ for communications satellites supporting cislunar operations? The present work seeks to address this question. 

The simplest model introducing the periodic forcing of the Sun in the Earth-Moon system is the bicircular restricted 4-body problem (BCP). This model is formed by assuming that the Earth-Moon barycenter travels around the Sun in a circular orbit at a constant rate, while the Earth and Moon move around their barycenter in a circular orbit. The BCP accounts for the Sun's gravitational effect on the particle (the direct effect) but neglects the gravitational effect on the relative motion of the Earth and Moon (the indirect effect). Due to this, the Sun, Earth, and Moon do not form a solution to the 3-body problem, hence we say the model is incoherent. Despite its incoherence, the BCP captures the qualitative behavior of the periodic solutions around $L_4$ similarly to its coherent counterpart, the quasi-bicircular restricted 4-body problem (QBCP), as shown in Figure \ref{fig:BCP_QBCP_L4}, borrowed from Andreu \cite{andreu1998quasi}. There has been some investigation into the dynamics around the triangular points in the BCP, namely the continuation of the $L_4$ point into its periodic counterparts \cite{simo1995bicircular}, as well as the computation of quasi-periodic invariant 2-tori in the vertical direction \cite{castella2000vertical}. Moreover, Jorba investigated the existence of stable motions near the triangular points of the Earth-Moon system in a high-fidelity ephemeris model, relating to work done in the BCP \cite{jorba2000numerical}. By a similar group of researchers, there has been more investigation modeling the dynamics affecting the Kordylewski clouds \cite{gimeno2024effect}. Along the lines of communications satellites being placed near $L_4$, some authors have used the BCP as their dynamical model for transfer design \cite{tan2020single,liang2021leveraging}.

In this work we use the Hill restricted 4-body problem (HR4BP) to model the motion of a particle in the Sun-Earth-Moon (SEM) system. The HR4BP is a coherent periodic Hamiltonian system that is both a generalization of the Earth-Moon CR3BP and Sun-Earth Hill's problem \cite{scheeres1998restricted,peterson2023vicinity}. The Sun-Earth-Moon HR4BP is coherent and models both the direct and indirect effect of the Sun, using a Hill's variational orbit as the relative motion of the Earth and Moon to capture the indirect effect. The HR4BP has been studied minimally in comparison to the BCP and QBCP \cite{scheeres1998restricted,olikara2016note,olikara_scheeres_2017,peterson2023vicinity,henry2023quasi}, and the purpose of employing this model in the present work is twofold. On one hand, only the original Scheeres paper which introduced the HR4BP has an investigation of the dynamics around EM $L_4$ in this model, and even that analysis is austere. So, we would like to further our understanding of the dynamics in this region of the phase space of the system. On the other hand, we would like to compare our results where applicable with the work previously done in the BCP and extend computations to include more than what has been done in the BCP to date. Additionally, from an astrodynamics perspective, transport between Earth and the Moon is of particular interest, so finding any avenues that exist passing through the triangular points may have value for transferring a satellite between cislunar regions of interest. Furthermore, practical stability regions near the triangular points have positive implications for the observability of spacecraft in cislunar space.

The paper is organized as follows. In Section \ref{sec:DynamicalModel} we review the dynamical model used in our investigation, including the continuation of equilibria and resonant periodic orbits near EM $L_4$. In Section \ref{sec:Methods} we give high-level overviews of the several techniques employed to study dynamics around the triangular points. Section \ref{sec:Results} makes up the majority of the text, as we detail our findings in the form of families of invariant tori, a Poincar\'e section, and transport via hyperbolic invariant manifolds. Finally, we draw comparisons where applicable between the existing work in the BCP and our examination of the HR4BP in Section \ref{sec:Comparison}.

\section{Dynamical Model}\label{sec:DynamicalModel}

\subsection{Hill Restricted 4-Body Problem}\label{sec:HR4BP}
The HR4BP is a coherent non-autonomous $\pi$-periodic generalization of the CR3BP and the Hill restricted 3-body problem. The derivation of the HR4BP is fully detailed in \cite{scheeres1998restricted}. There are three primary assumptions made in the HR4BP: $M_0 \gg M_1, M_2$, where $M_i$ is the mass of each primary body; the mass of the particle is negligible, i.e., $M_3 \ll 1$; and the distance between $M_1$ and $M_2$ is much less than their distance to $M_0$. In the SEM system, these assumptions are valid, where $M_0$ is the Sun, $M_1$ the Earth, and $M_2$ the Moon. Under these assumptions, the Sun and Earth-Moon barycenter move in a circular orbit and the relative motion of the Earth and Moon follows a specific solution to the Hill 3-body problem, a member of a 1-parameter family of planar periodic orbits known as Hill's variational orbits.

These periodic orbits depend on a parameter, $m$, which is the synodic period of the Earth-Moon orbit in years. Note for the SEM system, we take $m$ = 0.0808. In this work, we take the reference frame to be the average co-rotating frame of the smaller two primaries, rotating at a constant rate: $1 + 1/m$, where $1$ is the normalized rotation rate of the Sun and Earth-Moon barycenter frame, and $1/m$ is the average rotation rate of the Earth-Moon orbit. In addition to the parameter $m$, the HR4BP depends on a second parameter, $\mu$, the mass ratio of the smaller two primaries. Additionally, we take the following SEM unit normalizations. We normalize the length in the Earth-Moon frame by the mean distance between the Earth and the Moon (384,400 km). We normalize time $\tau$ such that $2\pi$ time units correspond to a synodic month (29.53 days), similar to \cite{olikara_scheeres_2017}.

To write the equations of motion, we define the conjugate momenta as $p_x = \dot{x} - (1+m)y$, $p_y = \dot{y} + (1+m)x$, and $p_z = \dot{z}$. The motion of the infinitesimal particle in the Earth-Moon rotating frame is then described as Hamiltonian dynamical system given by:
\begin{equation}\label{eqn:H_HR4BP}
    H = \frac{1}{2}(p_x^2 + p_y^2 + p_z^2) + \frac{1}{2}(1+m)^2(x^2+y^2) + (1+m)(yp_x-xp_y) - V(x,y,z,\tau;\mu,m),
\end{equation}
where $V(x,y,z,\tau;\mu,m)$ is given by:
\begin{align*}
V(x,y,z,\tau;\mu,m) &= \frac{1}{2}\left( 1 + 2m + \frac{3}{2}m^2 \right) (x^2 + y^2) - \frac{1}{2}m^2z^2 \\ &+ \frac{3}{4}m^2\big((x^2-y^2)\cos 2\tau - 2xy \sin 2\tau \big) \\
&+ \frac{m^2}{a_0^3(m)}\left( \frac{1-\mu}{r_{1-\mu}} + \frac{\mu}{r_\mu} \right),
\end{align*}
with $r_{1-\mu}$ and $r_\mu$ given by:
\begin{align*}
    r_{1-\mu}^2 &= \big( x + \mu(1+\bar{\xi}) \big)^2 + \big( y + \mu\bar{\eta} \big)^2 + z^2, \\
    r_{\mu}^2 &= \big( x - (1-\mu)(1+\bar{\xi}) \big)^2 + \big( y - (1-\mu)\bar{\eta} \big)^2 + z^2,
\end{align*}
where $\bar{\xi}$ and $\bar{\eta}$ are the normalized Hill variation orbit \cite{wintner1941analytical}:
\begin{align*}
    \bar{\xi}(\tau;m) &= \sum_{n = 1}^\infty \left( \frac{a_n(m)}{a_0(m)} + \frac{a_{-n}(m)}{a_0(m)} \right) \cos 2n\tau, \\
    \bar{\eta}(\tau;m) &= \sum_{n = 1}^\infty \left( \frac{a_n(m)}{a_0(m)} - \frac{a_{-n}(m)}{a_0(m)} \right) \sin 2n\tau,
\end{align*}
where the coefficients $a_n(m)$ can be computed following \cite{wintner1941analytical}; moreover, we use the coefficients given in \cite{olikara_scheeres_2017}. Note that the equations of motion are $\pi$-periodic in $\tau$, which corresponds to synodic half-month periodicity. Recall that $T = \frac{2\pi}{\omega_S}$, so, as the model is $\pi$-periodic, we say that $\omega_S = 2$. As $a_0(m) = m^{2/3}\left(1 - \frac{2}{3}m + \mathcal{O}(m^2) \right)$, we note that if $m \to 0$, the Hamiltonian of the restricted three-body problem is recovered.

We will use both the Hamiltonian and Newtonian formulations in this work, so we present the equations of motion in a Newtonian way:
\begin{equation}
    \ddot{x} - 2(1+m)\dot{y} = V_x, \quad \ddot{y} + 2(1+m)\dot{x} = V_y, \quad \ddot{z} = V_z,
\end{equation}
where $V$ is as defined in the Hamiltonian formulation above.

Finally, we mention several symmetries present in this model. Using the Hill approximation, the perturbing effect of the Sun is symmetric across the Earth-Moon barycenter, so the equations of motion are invariant under the transformation $\tau \mapsto \tau + \pi$. Additionally, the Earth-Moon orbit is $\pi$-periodic in this frame, so the system is $\pi$-periodic. The equations of motion are unchanged under the transformations:
\begin{align}
    S_y: (x,y,z,\dot{x},\dot{y},\dot{z},\tau) &\mapsto (x,-y,z,-\dot{x},\dot{y},-\dot{z},-\tau) \\ 
    S_z: (x,y,z,\dot{x},\dot{y},\dot{z},\tau) &\mapsto (x,y,-z,\dot{x},\dot{y},-\dot{z},\tau).
\end{align}
Thus, it suffices to consider the dynamics around $L_4$, as the dynamics around $L_5$ are obtained by applying the $S_y$ symmetry. Observe that this similar symmetry similarly holds in the CR3BP.

\subsection{Continuation of Equilibria and Resonant Periodic Orbits}\label{sec:Continuation}
The lowest dimensional invariant manifolds in the CR3BP are equilibrium points. We introduce the $\pi$-periodic forcing of the Sun by incrementally increasing $m$ from $m = 0$ to $m = 0.0808$, i.e., from the EM CR3BP to the SEM HR4BP. Once the periodic forcing is added, the equilibria become $\pi$-periodic orbits in general. Of course, there may be bifurcations of this periodic orbit in the continuation process, as seen around $L_2$ in the HR4BP \cite{olikara2016note,henry2023quasi}. We initialize the continuation of the $L_4$ point into its $\pi$-periodic orbit replacement taking $x_0 = x_{L_4}$, i.e., the initial state at $\tau = 0$ is taken as the CR3BP $L_4$ point. Using a standard single shooting algorithm provides the differential corrected periodic orbit \cite{seydel_practical_2010}.

Periodic orbits replacing the CR3BP equilibria are not the only periodic solutions of the periodically-forced system, as CR3BP periodic orbits of a resonant period may persist. Whereas the initial point at $\tau = 0$ is taken as the equilibrium point, we must carefully select the initial point for the continuation of a CR3BP periodic orbit into the HR4BP. In the SEM HR4BP, there is a 2:1 resonant periodic orbit which can be continued from the EM CR3BP. The 2:1 resonant orbit comes from the $2\pi$-periodic orbit of the $L_4$ short-period family in the CR3BP. We utilize sub-harmonic Melnikov theory to predict the points that persist \cite{guckenheimer2013nonlinear}. We utilize the energy principle--the change in energy from the start and end of a periodic orbit must be zero--following Cenedese and Haller \cite{cenedese2020conservative}, to construct a Melnikov function necessary for the subharmonic analysis at hand. 

To introduce the Melnikov used in this analysis, we write the HR4BP in the Newtonian formulation as $\ddot{\bm{r}} = \bm{f}(\bm{X}) + \varepsilon \bm{g}(\bm{X},\tau)$, where $\bm{f}(\bm{X})$ is the unperturbed system, $\bm{g}(\bm{X},\tau)$ is the $\pi$-periodic perturbation, and $\varepsilon > 0$ is small. Expanding the HR4BP equations of motion about $m = 0$, we have $\varepsilon \bm{g} = m\bm{g}_1 + m^2\bm{g}_2 + m^3\bm{g}_3 + \cdots$. Brown et al. articulated that $\bm{g}_1 \equiv 0$, so we use $\bm{g}_2$ as the next lowest order perturbation \cite{brown_2024_melnikov}. Repeating the computation of Brown et al., and using $\alpha \coloneqq \tau_0 + \tau$, we have:
\begin{align}
    \bm{g}_2 &= 3\bm{\Omega} \times \bm{r}' + \frac{3}{2} \bm{f} + \bm{h}_2, \\
    \bm{h}_2(\bm{X},\alpha) &= -\frac{3}{2}\begin{bmatrix}
        -\cos(2\alpha) & \sin(2\alpha) & 0 \\ \sin(2\alpha) & \cos(2\alpha) & 0 \\ 0 & 0 & -2/3
    \end{bmatrix} \bm{r} - \frac{1}{8}(1-\mu)\mu \begin{bmatrix} P \end{bmatrix} \begin{bmatrix}
        -8\cos(2\alpha) \\ 11\sin(2\alpha) \\ 0
    \end{bmatrix}, \\
    \begin{bmatrix} P \end{bmatrix} &= \left( \frac{1-\mu}{R_{1-\mu,C}^3} - \frac{\mu}{R_{\mu,C}^3} \right) [I_{3\times 3}] - \frac{3}{R_{1-\mu,C}^5} \bm{R}_{1-\mu,C} \bm{R}_{1-\mu,C}^T + \frac{3}{R_{\mu,C}^5} \bm{R}_{\mu,C} \bm{R}_{\mu,C}^T,
\end{align}
where $\mathbf{R}_{1-\mu,C}$ and $\mathbf{R}_{\mu,C}$ are the relative positions of the particle with respect to the Earth and the Moon, respectively, in the CR3BP. Originally derived in Brown et al. \cite{brown_2024_melnikov}, the Melnikov function for the 2:1 resonant periodic orbit in the HR4BP then takes the form:
\begin{equation}
    \mathcal{M}(s,\tau_0) = \int_0^{2\pi} \bm{h}_2( \bm{X}(s+\tau;\varepsilon = 0), \tau_0 + \tau) \cdot \bm{r}'(s+\tau;\varepsilon = 0) \, \text{d}\tau,
\end{equation}
where $\bm{X}(s+\tau;\varepsilon = 0)$ is the state along the orbit (where $s$ parameterizes the initial condition at $\tau_0$) in the unperturbed problem. The simple zeros of $\mathcal{M}$ represent the initial points along the orbit where no work is done by non-conservative forces over the orbit--a criterion for periodicity in the non-autonomous system. The top plot in Figure \ref{fig:L4_Continuation_Melnikov} shows the Melnikov function evaluated at each point along the 2:1 resonant orbit in the unperturbed problem. As there are 4 zeros, we continue the periodic orbit into the HR4BP from these points--observe that they are evenly spaced along the orbit and that opposite points correspond to the same periodic orbit with a phase shift. The bottom plot in Figure \ref{fig:L4_Continuation_Melnikov} shows the continuation in $m$ from the EM CR3BP ($m = 0$) into the SEM HR4BP ($m = 0.0808$). Notice that there is a period-doubling bifurcation between Branches A and B which destroys Branch B. Hence, only Branches A and C persist in the SEM system. Figure \ref{fig:L4_PeriodicOrbits} shows the surviving periodic orbits in the configuration space. We call Branch A the ``dynamical equivalent of $L_4$'' or DE $L_4$ for short; we call Branch C the 2:1 resonant orbit. 


The periodic forcing of the Sun qualitatively changes the linear normal behavior of $L_4$. So, we present in Table \ref{tab:L4_Stability} the stability data of the periodic orbits in the SEM HR4BP, i.e., the eigenvalues of the monodromy matrices, as well as the stability of the $L_4$ CR3BP equilibrium point and 2$\pi$-periodic orbit. The $\pi$-stroboscopic map is symplectic so the monodromy matrices are symplectic, hence the eigenvalues come in reciprocal pairs. Observe that the $L_4$ equilibrium point in the EM CR3BP is elliptic. In the SEM HR4BP, the elliptic $L_4$ point is replaced by a partially hyperbolic DE $L_4$ $\pi$-periodic orbit, i.e., having stability type center $\times$ saddle $\times$ center. This means there is a 1-parameter family of 2-dimensional invariant tori in the planar and vertical directions, and by coupling oscillations, there is a 2-parameter family of 3-dimensional tori. Conversely, the $2\pi$-periodic orbit in the EM CR3BP is elliptic and continues into an elliptic periodic orbit in the SEM HR4BP, though picks up a center mode because the unity eigenvalues transition into an additional center mode due to the periodic perturbation. This orbit--Branch C--is elliptic over the $m$ considered. Branch B, stemming from the same CR3BP periodic orbit, is partially hyperbolic throughout $m$ for which it persists. Branch B is killed in the bifurcation with Branch A, whence Branch A changes stability from elliptic to partially hyperbolic, inheriting the (un)stable manifold structure through the bifurcation. Figure \ref{fig:L4_Wsu_Evolution} shows the $\pi$-stroboscopic map of the (un)stable manifolds of Branches A and B for several values of $m$, namely before and after the bifurcation. In this figure, we see the structure of the (un)stable manifolds persists from Branch B to Branch A through the bifurcation. Note that these are not homoclinic points of the flow map, i.e., the (un)stable manifolds do not intersect in phase space. In the SEM HR4BP, Scheeres showed that the stable manifold of DE $L_4$ is confined to the interior of the structure, and the unstable manifold remains on the exterior.

\begin{figure}[H]
    \centering
    \includegraphics[width=0.77\linewidth]{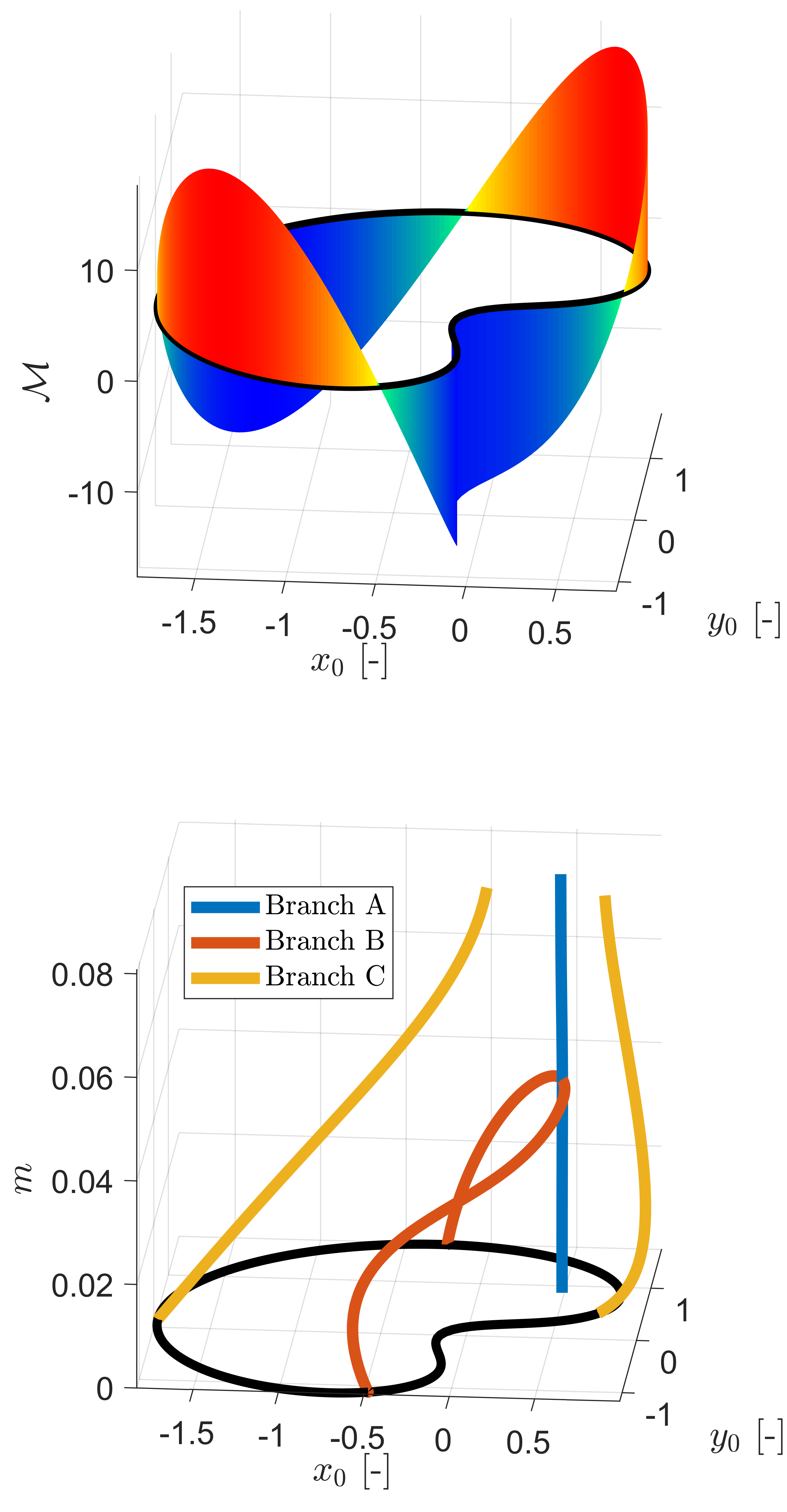}
    \caption{Continuation of 2:1 resonant periodic orbit, bifurcation with DE $L_4$ \cite{brown_2024_melnikov}.}
    \label{fig:L4_Continuation_Melnikov}
\end{figure}

\begin{figure}[H]
    \centering
    \includegraphics[width=0.8\linewidth]{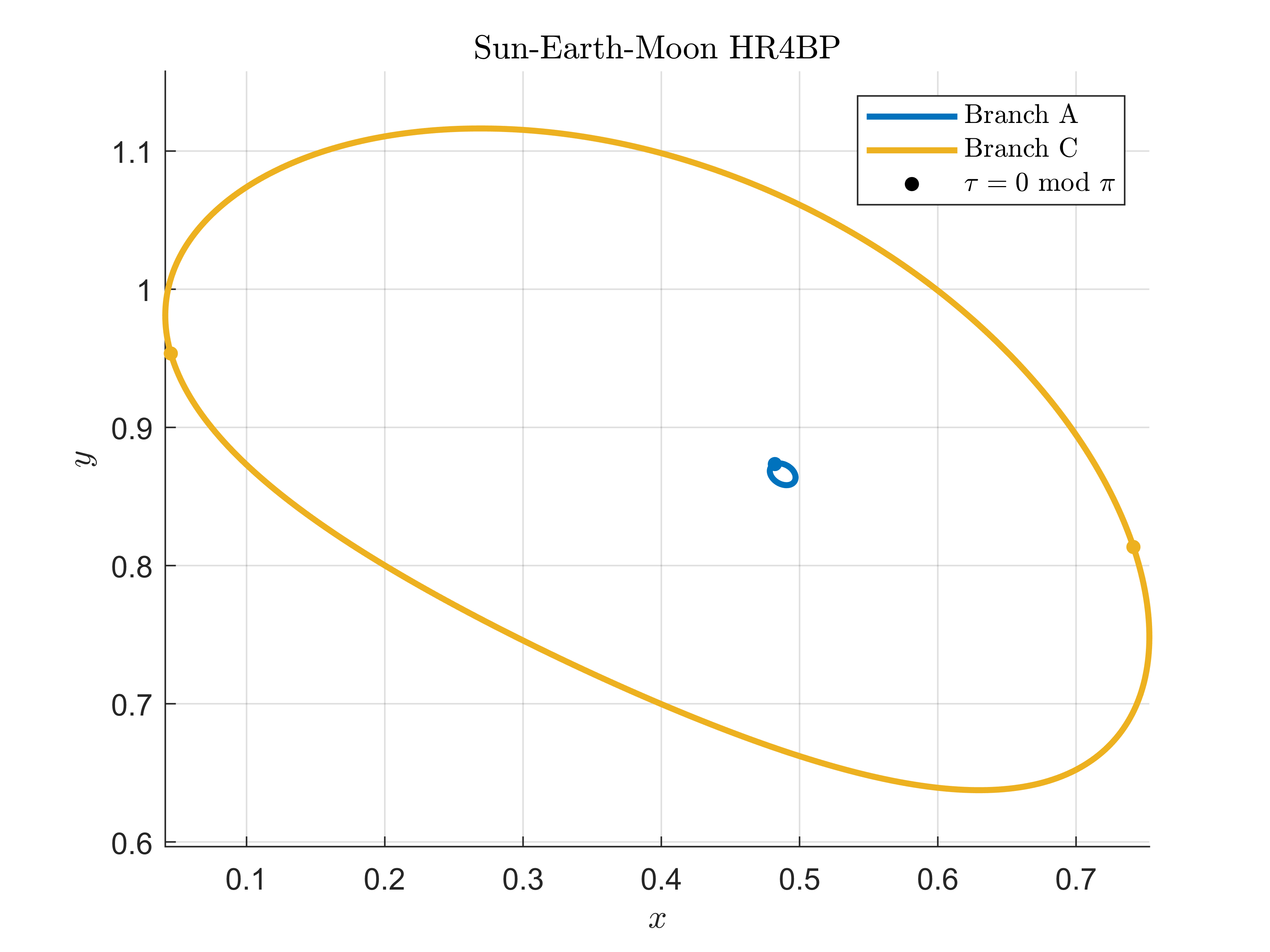}
    \caption{DE $L_4$ and 2:1 resonant periodic orbit in Sun-Earth-Moon HR4BP.}
    \label{fig:L4_PeriodicOrbits}
\end{figure}

\begin{figure}[H]
    \centering
    \includegraphics[width=0.9\linewidth]{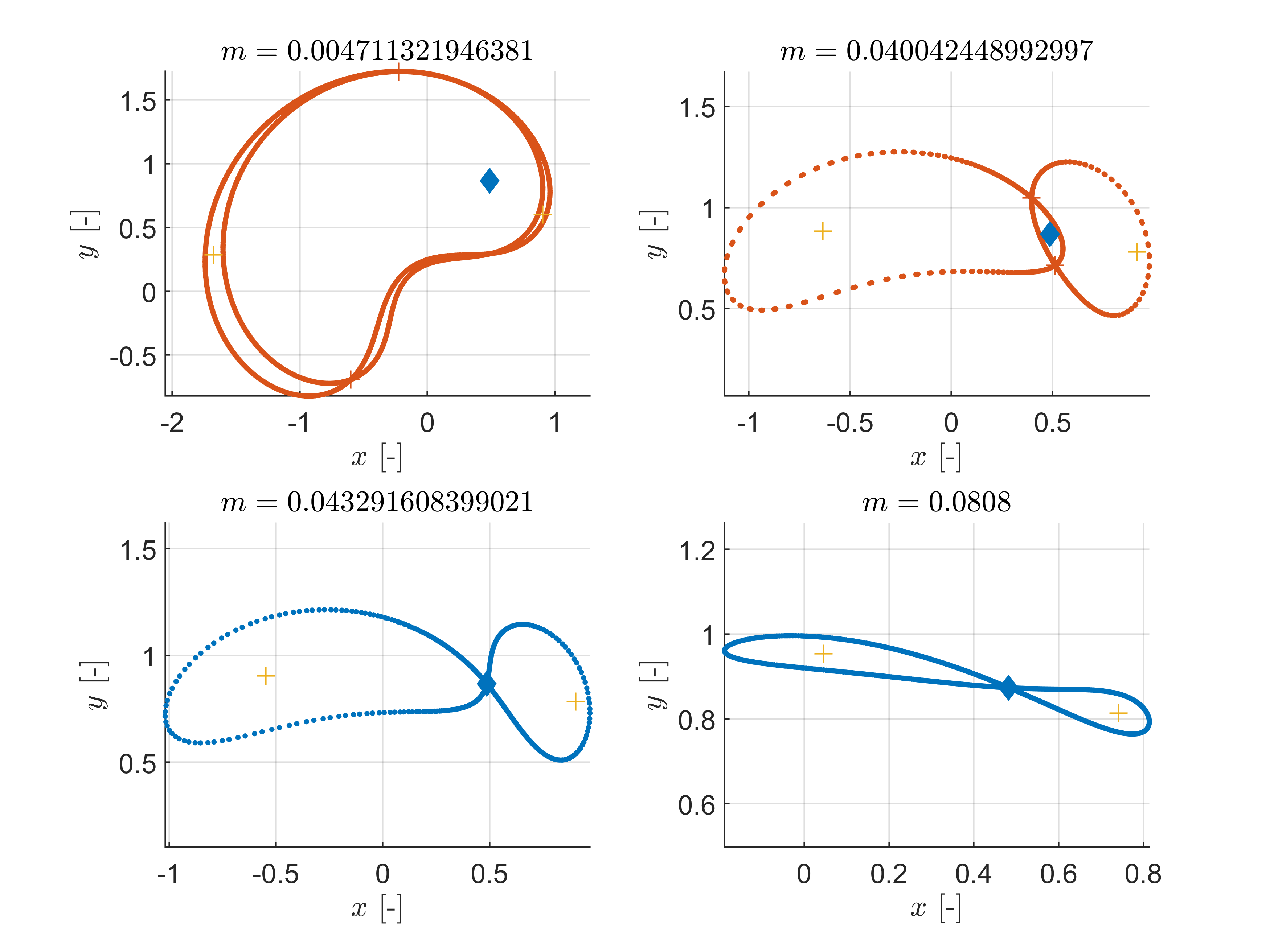}
    \caption{Stable and unstable manifolds of Branches B and A (DE $L_4$) under the $\pi$-stroboscopic map for several values of $m$. Note that the map shows $\tau = 0$ mod $\pi$.}
    \label{fig:L4_Wsu_Evolution}
\end{figure}

\begin{center}
\begin{table}[!htbp]
\caption{Stability of low-dimensional solutions in CR3BP and HR4BP}
\label{tab:L4_Stability}
\begin{tabular}{ccccc}
\hline
Solution & $\lambda_{1,2}$ & $\lambda_{3,4}$ & $\lambda_{5,6}$ \\
\hline
EM CR3BP $L_4$ EP & 0 $\pm \imag 0.9545$ & 0 $\pm \imag 0.2982$ & 0 $\pm \imag 1$ \\
EM CR3BP $2\pi$ PO & $0.6465 \pm \imag 0.7629$ & $1 \pm \imag 0$ & $0.9990 \pm \imag 0.0455$ \\
\hline 
SEM HR4BP DE $L_4$ & $0.5129 \pm \imag 0.8530$ & $-1.0547, -0.9482$ & $-0.9644 \pm \imag 0.2643$ \\
SEM HR4BP $L_4$ 2:1 & $-0.6043 \pm \imag 0.7968$ & $0.9885 \pm \imag 0.1513$ & $0.8534 \pm \imag 0.5212$ \\
\hline
\end{tabular}
\end{table}
\end{center}

\section{Methods}\label{sec:Methods}
In this section, we review the myriad of methods employed to study the dynamics near the Earth-Moon triangular points in the HR4BP. Firstly, we would like to characterize the bounded motions around DE $L_4$. Since DE $L_4$ is partially hyperbolic, a common strategy is to apply center manifold reduction wherein the bounded and hyperbolic motions are systematically decoupled up to some order of a series expansion around DE $L_4$. The primary limitation of the center manifold reduction is the radius of convergence which relates to the distance the expansion accurately represents the dynamics around DE $L_4$. Oftentimes center manifold reduction can be used in conjunction with the computation of quasi-periodic invariant tori to provide a more complete picture of the dynamics in a particular region; however, if the center manifold reduction technique is severely limited, this can give some insight into the dynamics qualitatively because there is a plausible explanation for this. Hence, we also apply the flow map method \cite{gomez_dynamics_2001,olikara2012numerical} to compute invariant 2-tori around DE $L_4$ and the 2:1 resonant orbit. We are also interested in the unbounded motions around DE $L_4$, so we compute the normal bundles of invariant tori. Using the stable and unstable bundles, we compute the global stable and unstable manifolds of particular 2-tori.

\subsection{Center Manifold Reduction}\label{sec:CMReduction}
The normal form computational techniques used here are based on \cite{jorba2020vicinity,rosales2023QBCP} and modified to the HR4BP, as in \cite{peterson2023vicinity}. We will give an overview of the fundamental aspects of the techniques; see the referenced texts for details on implementation. In this respect, the primary objective of this normal form technique is to re-center the Hamiltonian function about the dynamical equivalents of L1,2, expand as an infinite series, and, through successive time-dependent canonical transformations, change coordinates such that the system is autonomous and the local hyperbolic and elliptic modes of motion are decoupled. See the appendices of \cite{peterson2023vicinity} for the derivations used for re-centering and expanding the Hamiltonian around EM $L_{1,2}$--the procedure for EM $L_4$ is identical.

We re-center and expand the Hamiltonian in Equation \ref{eqn:H_HR4BP} about Earth-Moon $L_4$, denoted by $\mathcal{H}^{\text{EM}}$. See \cite{peterson2023vicinity} for details. This cancels terms of degree one in the Hamiltonian. In addition to re-centering, it is useful to autonomize the Hamiltonian by introducing a variable that is a symplectic conjugate to the periodic variable $\tau$, called $I_S$. Then, we define the autonomized Hamiltonian as:
\begin{equation}\label{eqn:H_auto}
    H(q,p,\tau,I_S) = \omega_S I_S + \mathcal{H}^{\text{EM}}.
\end{equation}
Note that the autonomization procedure above does not remove the dependence on $\tau$. We would like to remove the $\tau$-dependence up to a certain order, and this is reserved for a later step. 

Once the Hamiltonian function has been autonomized by introducing $I_S$ conjugate to $\tau$ and re-centered about the dynamical equivalent of interest to cancel terms of degree one, we utilize the Symplectic Floquet Theorem \cite{GomezJMS93rep} to put the quadratic terms into diagonal normal form with constant coefficients. Note that for degrees one and two, we need not expand the Hamiltonian as an infinite series of homogeneous polynomials \cite{andreu2002dynamics}. Applying the Floquet transformation to the Hamiltonian in Equation \ref{eqn:H_auto}, as well as the complexification given by
\begin{equation}\label{eqn:complexify}
    q_j = \frac{\tilde{q}_j + \imag \tilde{p}_j}{\sqrt{2}}, \quad p_j = \frac{\imag\tilde{q}_j + \tilde{p}_j}{\sqrt{2}}
\end{equation}
for $j = 2, 3$, we obtain the desired form:
\begin{equation}\label{eqn:H_cmplx}
    H(q,p,\tau,I_S) = \omega_S I_S + \omega_1 q_1 p_1 +  \imag\omega_2\tilde{q}_2\tilde{p}_2 + \imag\omega_3\tilde{q}_3\tilde{p}_3 + \sum_{n \geq 3} H_n(q,p,\tau,I_S).
\end{equation}

We now follow a similar discussion as \cite{jorba2020vicinity} to put the third-order terms into the desired form. The procedure is then applied similarly for subsequent orders. As mentioned above, we start with the Hamiltonian in Equation \ref{eqn:H_cmplx}, which is of the form: $H = \omega_S I_S + H_2 + H_3 + \cdots + H_n + \cdots$, where $H_2$ is in complex normal form and $H_n = H_n(q,p,\tau)$, $n \geq 3$, is a homogeneous polynomial of degree $n$ whose coefficients are complex-valued periodic functions of $\tau$. 

The normalization process begins at degree 3, using a generating function $G_3$ which is also a homogeneous polynomial of degree 3 in $(q,p)$ with coefficients depending periodically on $\tau$. Then, we continue at degree 4 with $G_4$, and so on. By choosing a generating function that defines a non-autonomous Hamiltonian vector field, and flowing along for unit time, we can define a time-dependent canonical transformation that removes the time dependence of the Hamiltonian in the new coordinates.

Starting with the following Hamiltonian:
\begin{equation}
    H = \omega_S I_S + H_2(q,p) + H_3(q,p,\tau) + \cdots + H_n(q,p,\tau) + \cdots,
\end{equation}
where $H_n(q,p,\tau) = \sum_{|k| = n} a_n^k(\tau) q^{k^1} p^{k^2}$, $a_n^k(\tau) = \sum_j a_{n,j}^k e^{j\imag\tau}$, and $k \in \mathbb{Z}^3 \times \mathbb{Z}^3$, denoted $k = (k^1,k^2)$, is a multi-index. Note that we use Fourier series to approximate the time-periodic coefficients $a_n^k(\tau)$. To remove certain terms from $H_3(q,p,\tau)$ (similarly for $H_n(q,p,\tau)$ for each $n$), we will make a change of variables generated by $G_3$:
\begin{equation}
    G_3 = G_3(q,p,\tau) = \sum_{|k|=3} g_3^k(\tau)q^{k^1}p^{k^2}, \quad g_3^k(\tau) = \sum_j g_{3,j}^k e^{j\imag\tau}.
\end{equation}
As $G_3$ generates a Hamiltonian vector field, the canonical transformation generated by $G_3$ is the time one flow along this vector field. We write explicitly the terms of degree 3 of the transformed Hamiltonian $\bar{H}_3$:
\begin{equation}\label{eqn:H3}
    \bar{H}_3 = H_3 - \omega_S \frac{\partial G_3}{\partial \tau} + \sum_{|k|=3} \langle \bar{\omega}, k^2 - k^1 \rangle g_3^k(\tau) q^{k^1} p^{k^2},
\end{equation}
where we define the vector of coefficients of $H_2$ as $\bar{\omega} = (\omega_1, \imag\omega_2, \imag\omega_3)$. Note that we consider here the stability case of \textit{saddle} $\times$ \textit{center} $\times$ \textit{center}, but this procedure is easily generalized to other stability cases. Continuing, note that $\bar{H}_3$ and $G_3$ are unknowns in Equation \ref{eqn:H3}, and we solve for $G_3$ presupposing $\bar{H}_3$ to have a desired form. Imposing the condition that the form of $\bar{H}_3 = \sum_{|k|=3}h_3^k(\tau)q^{k^1}p^{k^2}$, and grouping all terms with the same $k$, we obtain the following set of linear differential equations with $g_{3}^k$ as unknowns to be selected later:
\begin{equation}
    \omega_S \frac{d g_{3}^k}{d\tau} - \langle \bar{\omega}, k^2 - k^1 \rangle g_{3}^k = a_{3}^k - h_{3}^k.
\end{equation}
We can solve the above differential equations explicitly as:
\begin{equation}\label{eqn:SolveDE}
    g_{3}^k(\tau) = \sum_{j \notin J_k} \frac{a_{3,j}^k - h_{3,j}^k}{j\imag\omega_S - \langle \bar{\omega}, k^2 - k^1 \rangle} e^{j\imag\tau}, \quad j \geq 0,
\end{equation}
where $J_k = \{ j \in \mathbb{Z} \mid j\imag\omega_S - \langle \bar{\omega}, k^2 - k^1 \rangle = 0 \}$ is the ``resonance module,'' i.e., the set of indices of resonant terms. In the Sun-Earth-Moon HR4BP, $J_k$ contains only the zero vector for each $k$. Note that if $j \in J_k$, we have to impose $a_{3,j}^k = h_{3,j}^k$ and are unable to remove the time dependence. But, for $j \notin J_k$, we can choose $h_{3,j}^k$ and hence choose $G_3$.

In Equation \ref{eqn:SolveDE}, the $h_{3,j}^k$ are determined for the application. We seek to study the center manifold by reducing to an autonomous system and constructing an integral for the saddle such that, if set to zero, we obtain the center manifold. This means that we should choose $h_{n,j}^k$ as follows.

First, to eliminate time dependence, we choose the values $h_{n,j}^k = 0$ for $j > 0$ ($j \notin J_k$). Notice that by additionally setting $h_{n,0}^k = 0$, we remove also the monomial associated with $a_{n,j}^k$. To construct the saddle integral, we remove all monomials with $k_1^1 \neq k_1^2$. This means that the transformed Hamiltonian (up to some order $N$) will have the form:
\begin{equation}
    \bar{H}(q,p,\tau) = H^N(q_1p_1,q_2,q_3,p_2,p_3) + \mathcal{R}(q,p,\tau),
\end{equation}
where $\mathcal{R}$ is the remainder containing homogeneous polynomials of degree greater than $N$. Note that $\bar{H}$ depends on the product $q_1p_1$, rather than $q_1$ and $p_1$ separately. By defining the action variable $I_1$ and a symplectic conjugate hyperbolic angle $\theta_1$--see the Appendix in \cite{peterson2023local} for investigation of $\theta_1$--we obtain a canonical change of variables:
\begin{equation}
    \bar{H}(q,p) = H^N(I_1,q_2,q_3,p_2,p_3) + \mathcal{R}(I_1,\theta_1,q_2,q_3,p_2,p_3,\tau).
\end{equation}
If we neglect the remainder, the truncated Hamiltonian $H^N$ does not depend on $\theta_1$. Hence, $I_1$ is a first integral of the system, called the ``saddle integral'' \cite{peterson2023local}. If we set $I_1 = 0$ (in particular, $q_1 = p_1 = 0$), then $H^N$ provides an $N^{\text{th}}$-order approximation of the center manifold. We will exploit the saddle integral further in the next section to study the normal hyperbolic behavior of the center manifold. 

Note that we could additionally remove monomials to construct two integrals for the center manifold--this is the Birkhoff normal form. The benefit of the Birkhoff normal form is that it is integrable; however, the downside is the decreased radius of convergence. We note the construction of the center manifold is not unique. Depending on which monomials survive and which monomials are removed, one can obtain a different representation, e.g., the center manifold of $L_2$ in the QBCP is computed differently in \cite{andreu2002dynamics} and \cite{rosales2023QBCP}.

Following this procedure, the variables of $H^N(0,q,p)$ will be complex variables in general. Hence, a canonical realification transformation is performed using the inverse of the complexification transformation given by Equation \ref{eqn:complexify}. Further, the nonlinear change of variables is obtained by the sequence of canonical transformations generated by $G_3, G_4, \ldots$. We perform this order-by-order and can transform each coordinate $q_j$ or $p_j$ via the transformation defined by
\begin{align}
    \tilde{q}_j &= q_j + \{q_j,G_k\} + \frac{1}{2!}\{\{q_j,G_k\},G_k\} + \cdots + \frac{1}{n!}L_{G_k}^n(q_j) + \cdots \\
    \tilde{p}_j &= p_j + \{p_j,G_k\} + \frac{1}{2!}\{\{p_j,G_k\},G_k\} + \cdots + \frac{1}{n!}L_{G_k}^n(p_j) + \cdots,
\end{align}
where $L_{G_k}^n$ is the $n$-th order Lie derivative. See \cite{jorba1999methodology} for details.

Finally, we investigate the dynamics in the local center manifold around DE $L_4$ by studying the vector field induced by $H^N(0,q,p): \mathbb{R}^4 \mapsto \mathbb{R}$ which is independent of time. In these dynamics, the energy integral is recovered. This reduces the problem of studying dynamics around DE $L_4$ to the study of area-preserving maps parameterized by the energy \cite{peterson2023vicinity}. Once the energy is fixed, we can take a particular Poincar\'e section to obtain a 2-dimensional picture of the local center manifold. In this work, we use the vertical section $\Sigma_v = \{ q_2 = 0 \}$, which corresponds to fixing $z = 0$ in synodical coordinates to first order \cite{jorba2020vicinity}.

\subsection{Computation of Invariant Tori}\label{sec:InvariantTori}
In this work, we compute families of 2-dimensional quasi-periodic invariant tori by computing a parameterization $v: \mathbb{T}^2 \to \mathcal{T} \subset \mathbb{R}^6$ that maps angles $\theta = (\theta_0,\theta_1)$ to a state on the quasi-periodic orbit, $\mathcal{T}$, that is invariant under the flow induced by the Hamiltonian function $H$. The dynamics on the standard 2-torus, $\mathbb{T}^2$, are linear and the constant frequencies of oscillation, $\omega = (\omega_0,\omega_1)$, are non-resonant. In particular, the frequency $\omega_0$ is the frequency of the perturbation, i.e., $\omega_0 = \omega_S = 2$. 

Since the parameterization $v$ is invariant under the flow, the invariant curve defined by the map $u(\theta_1) = v(0,\theta_1): \mathbb{T}^1 \mapsto \mathcal{T} \subset \mathbb{R}^6$, satisfies
\begin{equation}\label{eqn:invariance}
    \varphi_T(u(\theta_1)) = u(\theta_1 + \rho)
\end{equation}
where $\varphi_t$ is the flow induced by $H$ for time $t \in \mathbb{R}$, $T = 2\pi/\omega_0 = \pi$ is the period of the perturbation, and $\rho = 2\pi \omega_1/\omega_0$ is the rotation number of the torus. 

To compute $u$ explicitly, we discretize the invariant curve $u(\theta_1)$ over $2N+1$ angles whence the computation of invariant curves becomes a boundary value problem in which the $2N+1$ states must meet the boundary condition described by Equation \ref{eqn:invariance}. We solve the boundary value problem by defining the grid of angles, the grid of points along the invariant curve, and the free variable vector as
\begin{align}
    \theta_{1,k} &= \frac{2\pi k}{2N+1}, \quad k = 0,\ldots,2N, \\
    U &= \begin{bmatrix} u(\theta_{1,0})^{\intercal} & \cdots & u(\theta_{1,2N})^{\intercal} \end{bmatrix}^{\intercal}, \\
    z &= \begin{bmatrix} U^{\intercal} & \rho \end{bmatrix}^{\intercal},
\end{align}
respectively. Defining the constraint vector as
\begin{equation}
    g(z) \coloneqq R_{-\rho} \varphi_T(U) - U = 0,
\end{equation}
where we apply the rotation operator $R_{-\rho}$, sending arguments of functions depending on $\theta_1$ to $\theta_1 - \rho$, and the flow $\varphi_T$ is applied to each block element, as in \cite{henry2023quasi}, solutions to $g(z) = 0$ implicitly define a smooth manifold that is the 2-dimensional quasi-periodic invariant torus. It is important to note that the rotation operator can be constructed by composing discrete Fourier transform operations \cite{olikara2012numerical,henry2023quasi}. Note that we solve this boundary value problem using a multiple shooting scheme with a Newton update step. 

As any phase shift also satisfies the boundary value problem, we also add a phase constraint that minimizes the phase difference between the initial curve and the previously computed invariant curve, following \cite{olikara2012numerical,henry2023quasi}. To compute families of invariant 2-tori, which lie in Cantorian one-parameter families in the HR4BP \cite{jorba_normal_1997}, we also employ a pseudo-arclength continuation scheme \cite{seydel_practical_2010}. Hence, we add the continuation constraint
\begin{equation}
    c(z) \coloneqq \frac{1}{2N+1}(U - \tilde{U})^{\intercal}n_u + (\rho - \tilde{\rho})n_\rho - \Delta s = 0,
\end{equation}
where $\tilde{U}$ and $\tilde{\rho}$ denote the previously computed grid points along the invariant and rotation number, repsectively, $n_u$ and $n_\rho$ are approximations of the null space direction of the constraint Jacobian matrix corresponding to the invariant curve and rotation number, and $\Delta s$ is the step size.

\subsection{Normal Behavior of Invariant Curves}\label{sec:Stability}
We are interested in studying the stability of 2-dimensional quasi-periodic orbits in the HR4BP. By using the flow map method, we reduce the problem to that of computing the linearized normal behavior of an invariant curve of a diffeomorphism (the flow map $\varphi_T$) \cite{jorba_numerical_2001}. Given an invariant curve $u(\theta_1)$ such that $\varphi_T(u(\theta_1)) = u(\theta_1+\rho)$ for all $\theta_1 \in \mathbb{T}^1$. A small displacement $h \in \mathbb{R}^6$ with respect to a point $u(\theta_1)$ on the curve is:
\begin{equation}
    \varphi_T(u(\theta_1)+h) = \varphi_T(u(\theta_1)) + D_{u}\varphi_T(u(\theta_1))h + \mathcal{O}(\|h\|^2).
\end{equation}
The linear normal behavior is described by the system
\begin{equation}
    \bar{u} = A(\theta_1){u}, \quad \bar{\theta}_1 = \theta_1 + \rho,
\end{equation}
where $A(\theta_1) = D_{u}\varphi_T(u(\theta_1))$. 
The rotation operator defined in the previous section is applied to the generalized eigenvalue problem as $R_\rho: \psi(\theta) \mapsto \psi(\theta+\rho)$. We consider the generalized eigenvalue problem in which we look for pairs $(\lambda,\psi) \in \mathbb{C} \times (\mathcal{C}(\mathbb{T}^1,\mathbb{C}^n)\setminus\{0\})$ such that 
\begin{equation}\label{eqn:EigenvalueProblem}
    A(\theta_1)\psi(\theta_1) = \lambda R_\rho \psi(\theta_1).
\end{equation}
To compute the stability of invariant tori, we discretize the operator $R_{-\rho}\circ A(\theta_1)$, which is represented by a matrix. We then compute the eigenvalues and eigenvectors of this matrix by a standard numerical procedure. The resulting eigenvalues will form circles in the complex plane, covered by the same number of points as in the discretized invariant curve. If the invariant curve $u(\theta_1)$ is reducible, then its normal behavior will depend on six independent eigenvalues and eigenfunctions from this discretized set \cite{jorba_numerical_2001}. As the map $\varphi_T$ is symplectic, the corresponding 6 eigenvalues will come in reciprocal and conjugate pairs, and there will be one pair of unity eigenvalues. The eigenfunctions corresponding to the unity pair are the tangent vector of the invariant curve along the flow and the tangent vector pointing in the direction of the one-parameter family of invariant 2-tori. 

While the eigenvalues computed in this discretization will have different accuracy, due to the discretization, we implement a sorting method described in Jorba \cite{jorba_numerical_2001} and summarize the technique here. Representing each eigenfunction as a Fourier series $\sum_j \psi_j \text{exp}(\imag j\theta_1)$, we want to determine which eigenfunctions have the least discretization error. Hence, we consider the truncated Fourier series
\begin{equation}
    \psi_j^N(\theta_1) = \sum_{j = -N}^N \psi_j \text{exp}(\imag j\theta_1).
\end{equation}
The discretization error of the eigenfunctions is given by the magnitude of Fourier coefficients for $j > N$ \cite{henry2023quasi}. As the eigenfunctions are analytic, their Fourier coefficients should decay exponentially from where the spectrum is centered \cite{henry2023quasi}. So, by checking the truncation error defined by 
\begin{equation}
    TE(\psi,N) = \sum_{|j| > N} |\psi_j| |j|,
\end{equation}
we find the most ``centered'' eigenfunctions, and we choose these as our representative eigenfunctions and corresponding eigenvalues as the six representatives of the equivalence classes, provided the tails are sufficiently small. These representative eigenvalues and eigenfunctions then provide an accurate description of the normal dynamics. In particular, an eigenvalue pair on the unit circle in the complex plane describes a normally elliptic direction, and an eigenvalue pair on the real line away from $\pm 1$ describes a normally hyperbolic direction. Changes in stability type indicate a bifurcation in the family of quasi-periodic invariant tori, similar to Krein collisions of periodic orbits in Hamiltonian systems \cite{arnold_mathematical_2006}.

\subsection{Hyperbolic Invariant Manifolds of Invariant Tori}\label{sec:InvariantManifolds}
The stable manifold theorem gives the existence of stable and unstable invariant manifolds emanating from quasi-periodic invariant 2-tori having a pair of real eigenvalues with modulus different from 1. The (global) stable and unstable manifolds are defined as the collections of points that asymptotically approach the invariant 2-torus as $t \to +\infty$ and $t \to -\infty$, respectively. In the section above, we described a process for computing eigenvalues and eigenfunctions of invariant curves. We follow \cite{jorba2020transport,simo1990analytical} to compute the hyperbolic invariant manifolds of a normally hyperbolic invariant 2-torus. Suppose we have computed some invariant curve $u(\theta_1)$ which has a pair of hyperbolic eigenvalues $\lambda_{s,u} \in \mathbb{R}$ with corresponding eigenfunctions $\psi_{s,u}(\theta_1)$. For $h \in \mathbb{R}$ small and any $\theta \in [0,2\pi]$, these hyperbolic eigenvalues and eigenfunctions satisfy:
\begin{align}
    \varphi_T(u(\theta_1)+h\psi_{s,u}) &= \varphi_T(u(\theta_1)) + h A(\theta_1)\psi_{s,u}(\theta_1) + \mathcal{O}(\|h\|^2) \\
    &= u(\theta_1 + \rho) + h\lambda_{s,u}\psi_{s,u}(\theta_1 + \rho) + \mathcal{O}(\|h\|^2),
\end{align}
where $A(\theta_1)$ is defined in the previous section. Hence, the map
\begin{equation}
    (\theta_1,h) \mapsto u(\theta_1) + h \psi_{s,u}(\theta_1)
\end{equation}
is a parameterization of the linearization of the stable and unstable manifolds along the invariant curve $u(\theta_1)$. To compute the global stable and unstable manifolds of the invariant curve, one can fix a value of $h = h_0$ so that the error of the linearization is sufficiently small--we have used $h_0 = 10^{-5}$ so that the error is $\mathcal{O}(10^{-10})$. However, as the HR4BP is $\pi$-periodic, under the $\pi$-periodic flow map a point of the invariant curve perturbed by $h_0$ onto the stable (unstable, respectively) manifold will be perturbed by $h_0/\lambda_s$ ($\lambda_u h_0$, respectively). By choosing $h \in [h_0, h_0/\lambda_s]$ ($[h_0, \lambda_u h_0]$, respectively), we can effectively parameterize (in a linear sense) the stable and unstable manifolds of not only the invariant curve $u(\theta_1)$, but of the entire 2-torus $v(\theta_0,\theta_1)$ under the $\pi$-stroboscopic map. In other words, we have parameterized the fundamental cylinder of the hyperbolic invariant manifold. To globalize the manifold, we define a mesh of points along the cylinder $(\theta_1,h) \in [0,2\pi] \times [h_0, h_0/\lambda_s]$ ($[h_0,\lambda_u h_0]$, respectively) and propagate the points on the cylinder backward (forward, respectively) to compute the stable (unstable, respectively) manifold. Note that we perturb in $\pm h$ to obtain both sides of the invariant manifold. In this work, we choose a $251 \times 251$ grid for computations.

\section{Results}\label{sec:Results}
In this section, we display and analyze the results from the methods described in the previous section applied to the Sun-Earth-Moon HR4BP. Firstly, we include a Poincar\'e map generated via center manifold reduction. Limitations of this method applied to $L_4$ are discussed. Then, five families of Lyapunov invariant 2-tori are presented along with their normal behavior. Finally, we show the existence of natural trajectories that escape the region around the triangular points, escaping the system or encountering regions near the primaries.

\subsection{Poincar\'e Map}\label{sec:PoincareMap}
\begin{figure}
    \centering
    \includegraphics[width=0.9\linewidth]{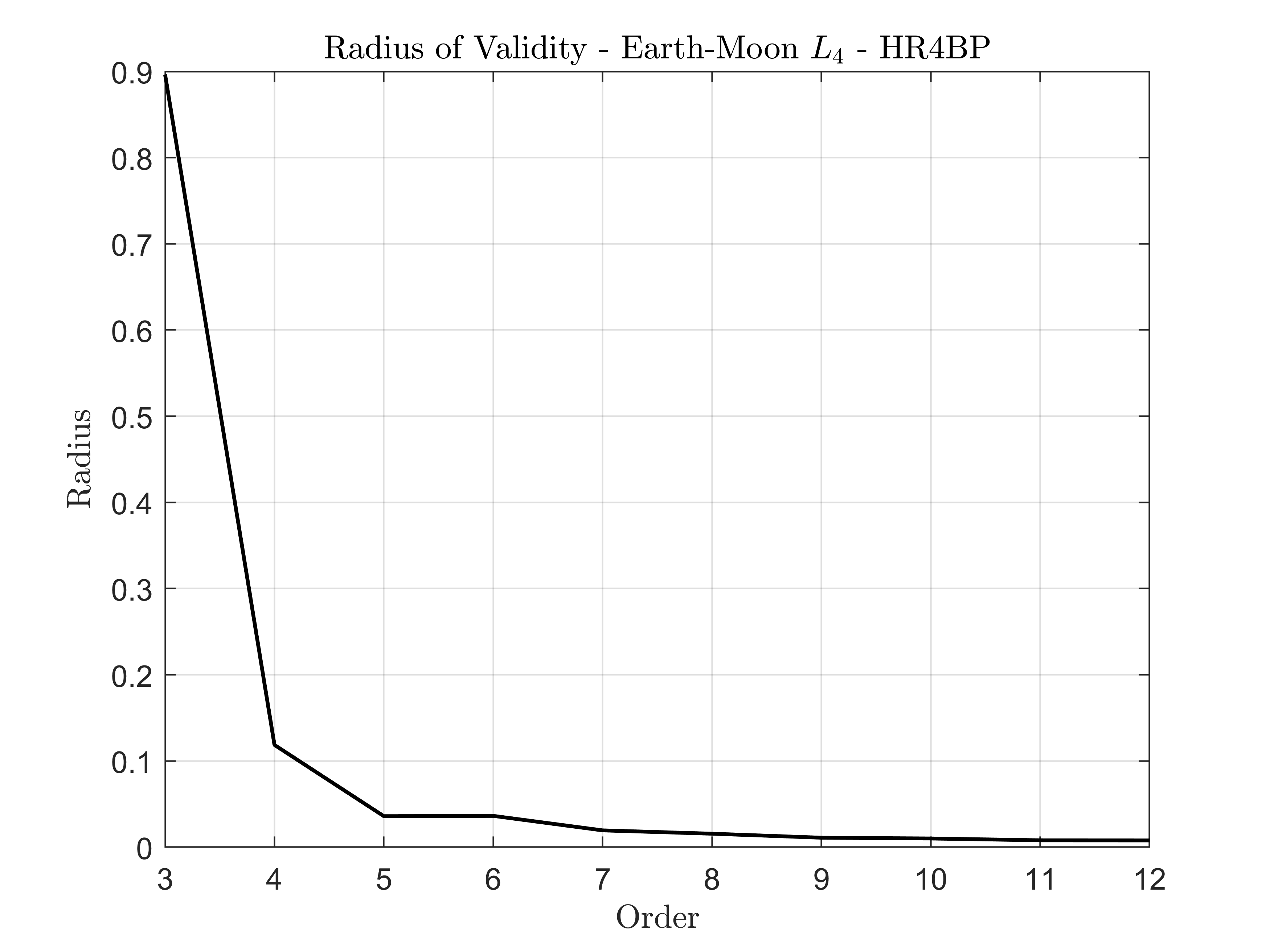}
    \caption{Radius of validity for center manifold reduction around DE $L_4$.}
    \label{fig:RoC}
\end{figure}

\begin{figure}
    \centering
    \includegraphics[width=0.4\linewidth]{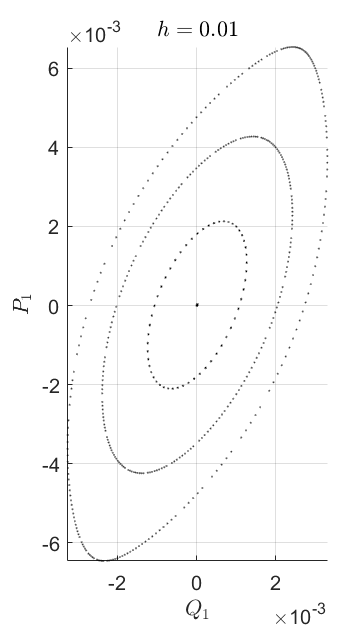}
    \caption{Poincar\'e section $\Sigma_v$ of degree 4 at fixed energy value $h = 0.01$.}
    \label{fig:PoincareMap}
\end{figure}

The center manifold reduction was computed about DE $L_4$ up to degree 12. Figure \ref{fig:RoC} shows the radius of validity for the center manifold reduction around DE $L_4$ up to degree 12. Observe that this radius decreases suddenly and asymptotically approaches zero. The small radius of convergence indicates that the Poincar\'e map will be severely limited in its utility to accurately describe dynamics near DE $L_4$. Namely, Figure \ref{fig:PoincareMap} shows a Poincar\'e map of the center manifold reduction of degree 4 at a fixed energy value $h = 0.01$. Note that the order and recovered energy integral are both small due to the divergence properties of the series expansion. Figure \ref{fig:PoincareMap} shows the existence of 2- and 3-dimensional planar Lyapunov invariant tori around DE $L_4$. Yet, the divergence properties of the series expansion limit the computation and numerical precision of the invariant tori. Hence, we utilize the flow map method to parameterize individual invariant 2-tori to study the invariant manifold structures around DE $L_4$ and $L_4$ 2:1 resonant periodic orbit.

\subsection{Families of Invariant 2-tori}\label{sec:Families}
In this section, we present the results of computing the five families of Lyapunov invariant 2-tori around DE $L_4$ and the $L_4$ 2:1 resonant periodic orbit and their normal behavior. We will identify bifurcations and compute bifurcations in the cases which are created by the perturbation. In other words, if the bifurcation exists between periodic orbits in the Earth-Moon CR3BP, then the bifurcation will be identified but not computed or continued. 

\subsubsection[DE L4 H]{DE $L_4$ H}
The first family of Lyapunov invariant 2-tori around DE $L_4$ is the planar family, which we denote as ``DE $L_4$ H.'' The ``H'' stands for horizontal, as opposed to ``V'' for vertical. By taking as an initial guess the excited planar center mode of DE $L_4$ periodic orbit, we compute this family of invariant 2-tori using pseudo-arclength continuation and the algorithm described in Section \ref{sec:Methods}. This family is qualitatively similar to (and a generalization of) the long-period family of periodic orbits around $L_4$ in the Earth-Moon CR3BP. 

\begin{figure}
     \centering
     \begin{subfigure}[b]{0.9\textwidth}
         \centering
         \includegraphics[width=\textwidth]{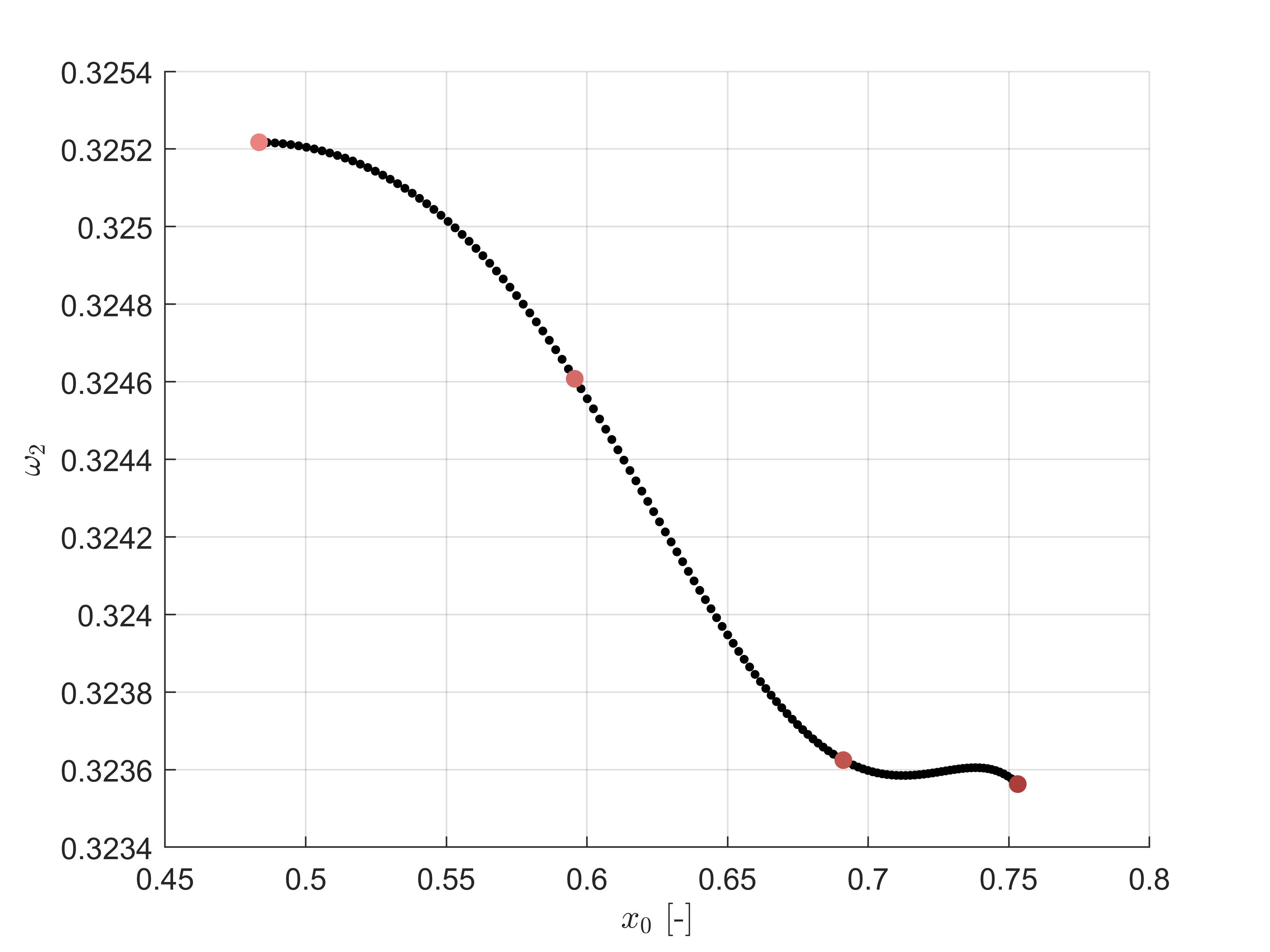}
     \end{subfigure}
     \hfill
     \begin{subfigure}[b]{0.9\textwidth}
         \centering
         \includegraphics[width=\textwidth]{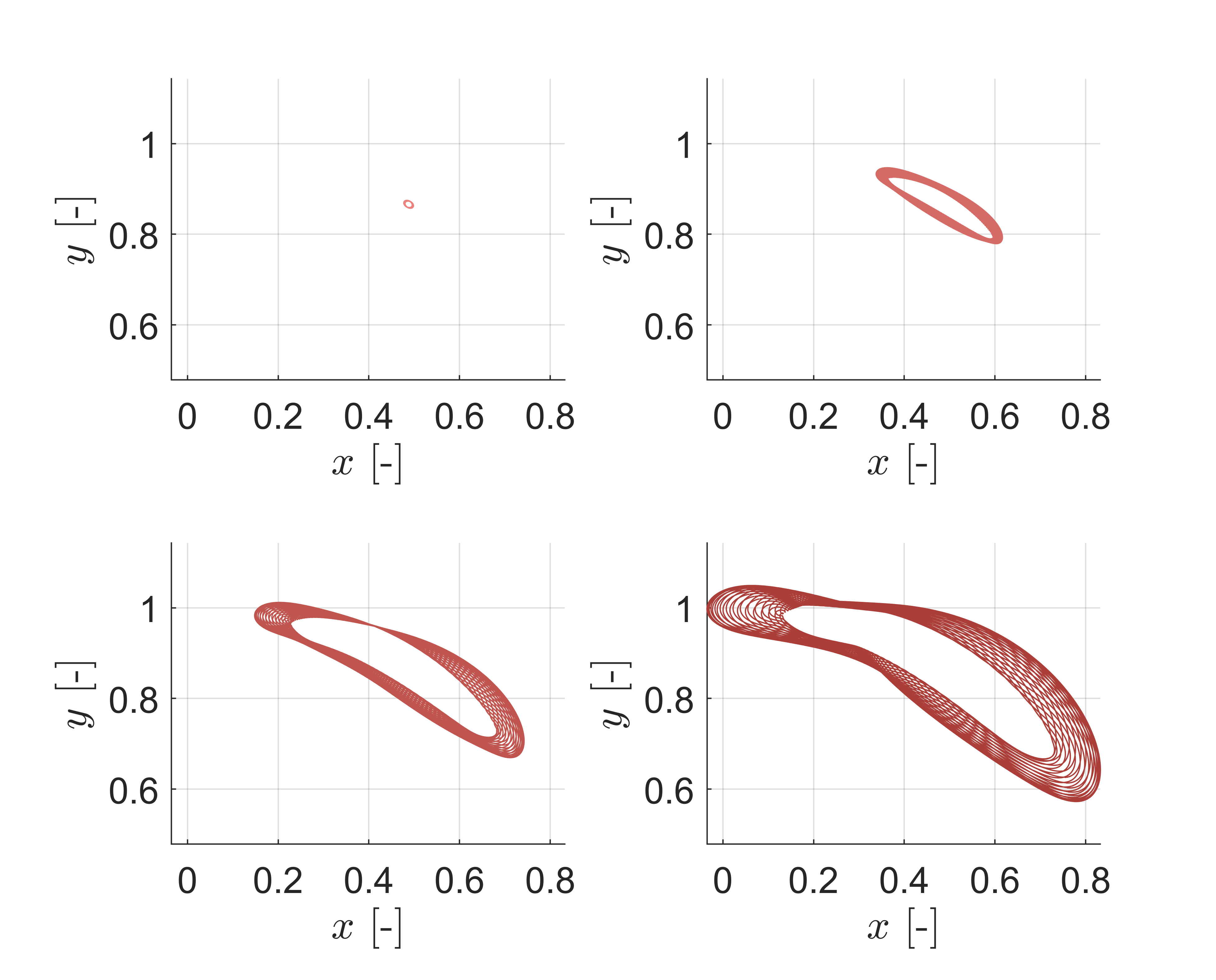}
     \end{subfigure}
        \caption{Hodograph and sample invariant tori of DE $L_4$ H family. See text for details.}
        \label{fig:DEL4_H}
\end{figure}

Figure \ref{fig:DEL4_H} shows a hodograph of the family computation, along with sample torus representations. The family grows as $x_0$ increases. We follow the convention to represent families of invariant tori as with periodic orbits, i.e., to choose a particular state variable--$x(\tau)$ at $\tau = 0 \mod \pi$--as a variable plotted against the frequency associated to the center mode. However, unlike in the case of periodic orbits, the choice of this variable is not unique--any point along the $\tau = 0 \mod \pi$ invariant curve is equivalent. We choose the first $x$-coordinate along this invariant curve given by the numerical output of the differential correction. This results in a smooth curve for the family's hodograph. The family computation ends when the algorithm fails to converge. In this instance, the $\tau = 0 \mod \pi$ invariant curve became less analytic, resulting in an increased number of Fourier modes used to describe the curve. While the family could be continued farther, we stopped the computation when the maximum number of Fourier modes was reached. Alternatively, one could use the invariant curve along the other angular variable. 

\begin{figure}
    \centering
    \includegraphics[width=\linewidth]{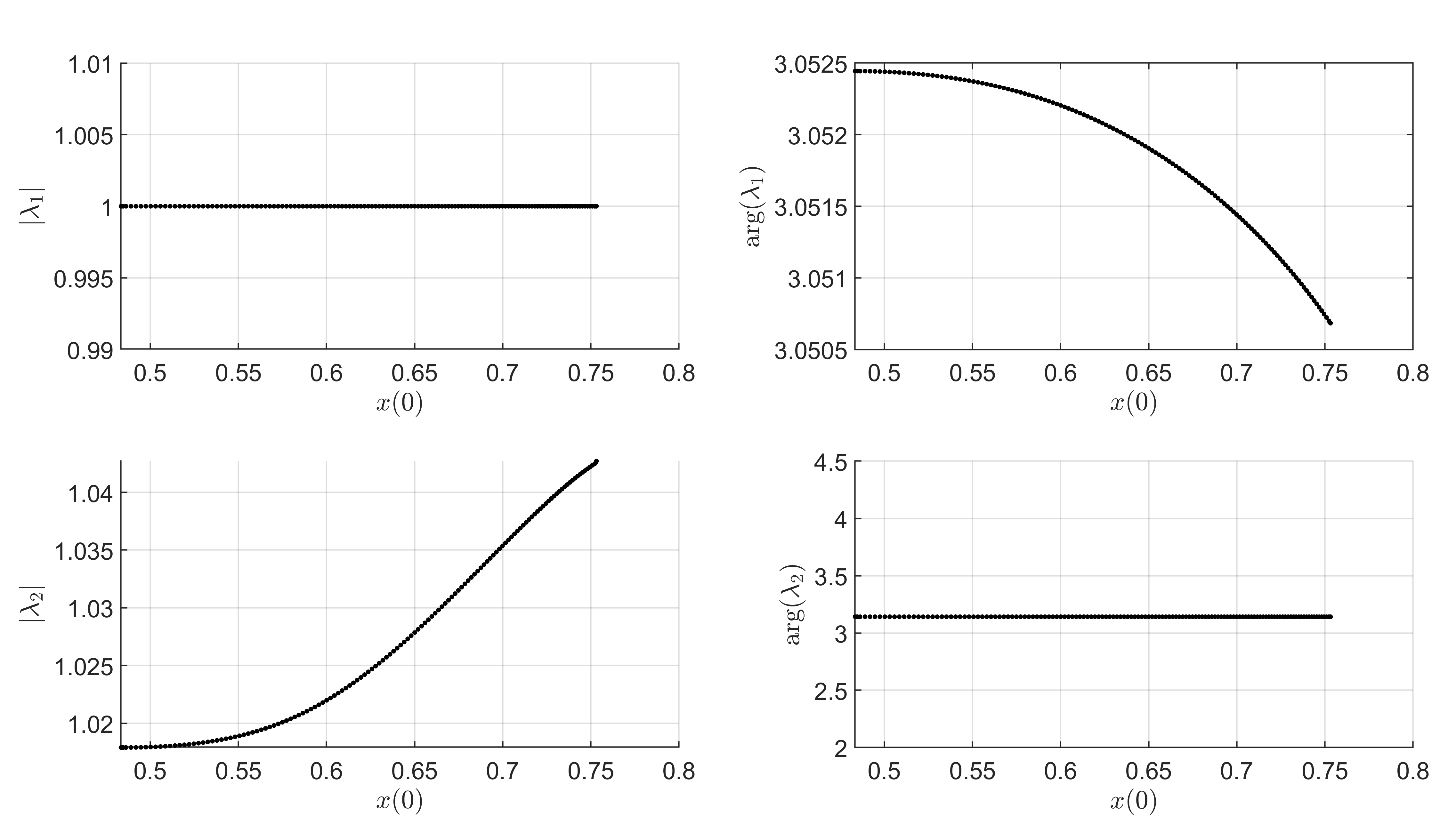}
    \caption{Normal behavior of DE $L_4$ H family of invariant tori.}
    \label{fig:Stab_DEL4_H}
\end{figure}

Figure \ref{fig:Stab_DEL4_H} shows the computed normal behavior of the family of tori. As in the previous figure, the family grows in increasing $x_0 = x(0 \mod \pi)$. To show the stability, we describe two eigenvalues representing each normal mode. Due to the symplectic nature of the system, these eigenvalues come in reciprocal pairs (Krein quartets, in fact), so we pick one from each pair. The normal behavior for the DE $L_4$ H family remains of type center $\times$ saddle, as seen most clearly by the magnitude of eigenvalues $\lambda_1$ and $\lambda_2$, the center and unstable eigenvalues, respectively. The behavior of this family is the most straightforward of the five presented, and we will return to DE $L_4$ H when considering transport near $L_4$ in the HR4BP in a later section. 

\subsubsection[DE L4 V]{DE $L_4$ V}
The DE $L_4$ V family is a collection of invariant 2-tori emanating from the center mode of DE $L_4$ in the vertical ($z$) direction. This family is qualitatively similar to (and a generalization of) the vertical Lyapunov family of periodic orbits around $L_4$ in the Earth-Moon CR3BP. Using the same computational methods as above, we computed invariant 2-tori until a maximum $z$-value was reached because we are interested in the dynamics within this bound. 

\begin{figure}
     \centering
     \begin{subfigure}[b]{0.8\textwidth}
         \centering
         \includegraphics[width=\textwidth]{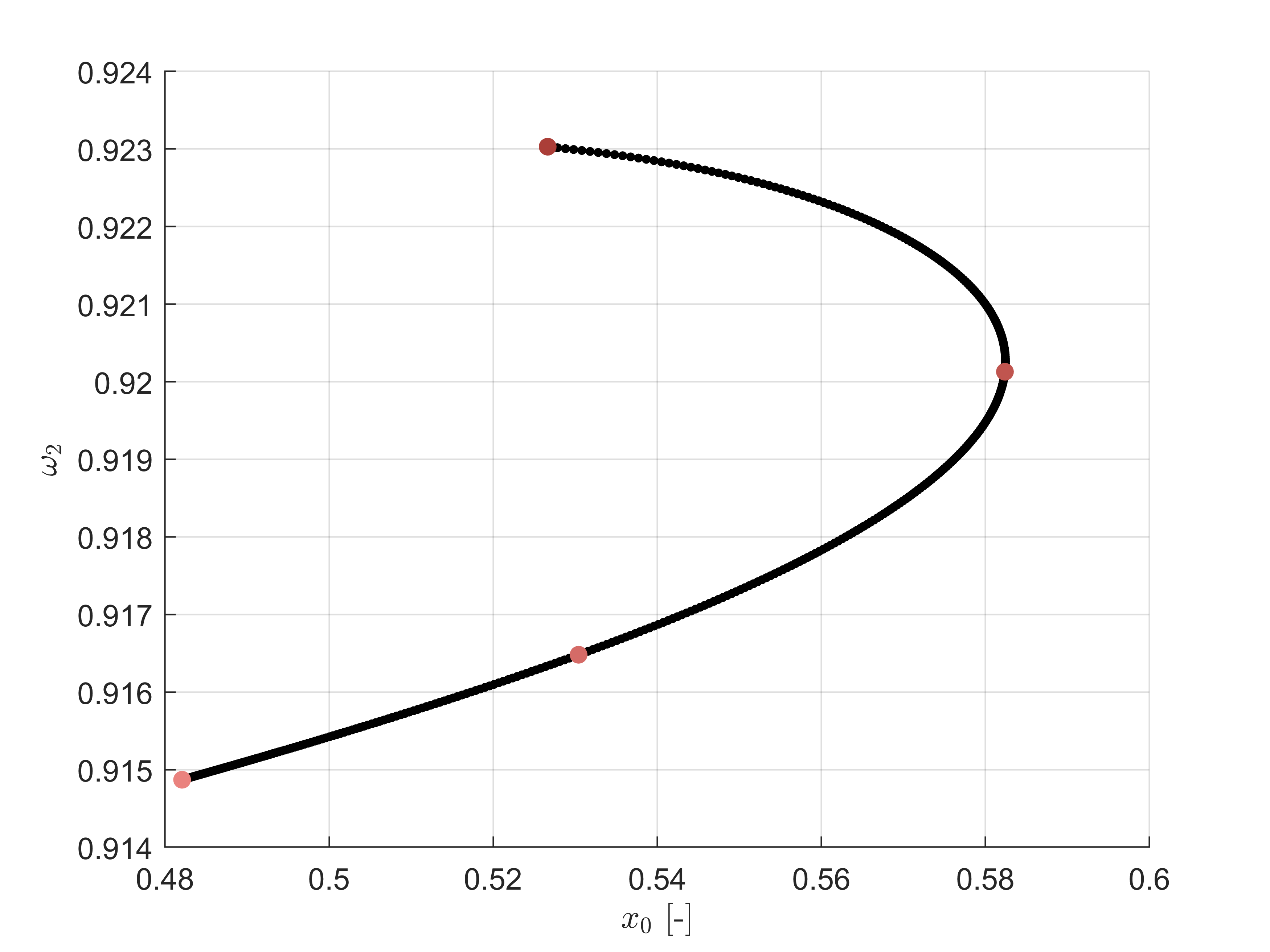}
     \end{subfigure}
     \hfill
     \begin{subfigure}[b]{0.8\textwidth}
         \centering
         \includegraphics[width=\textwidth]{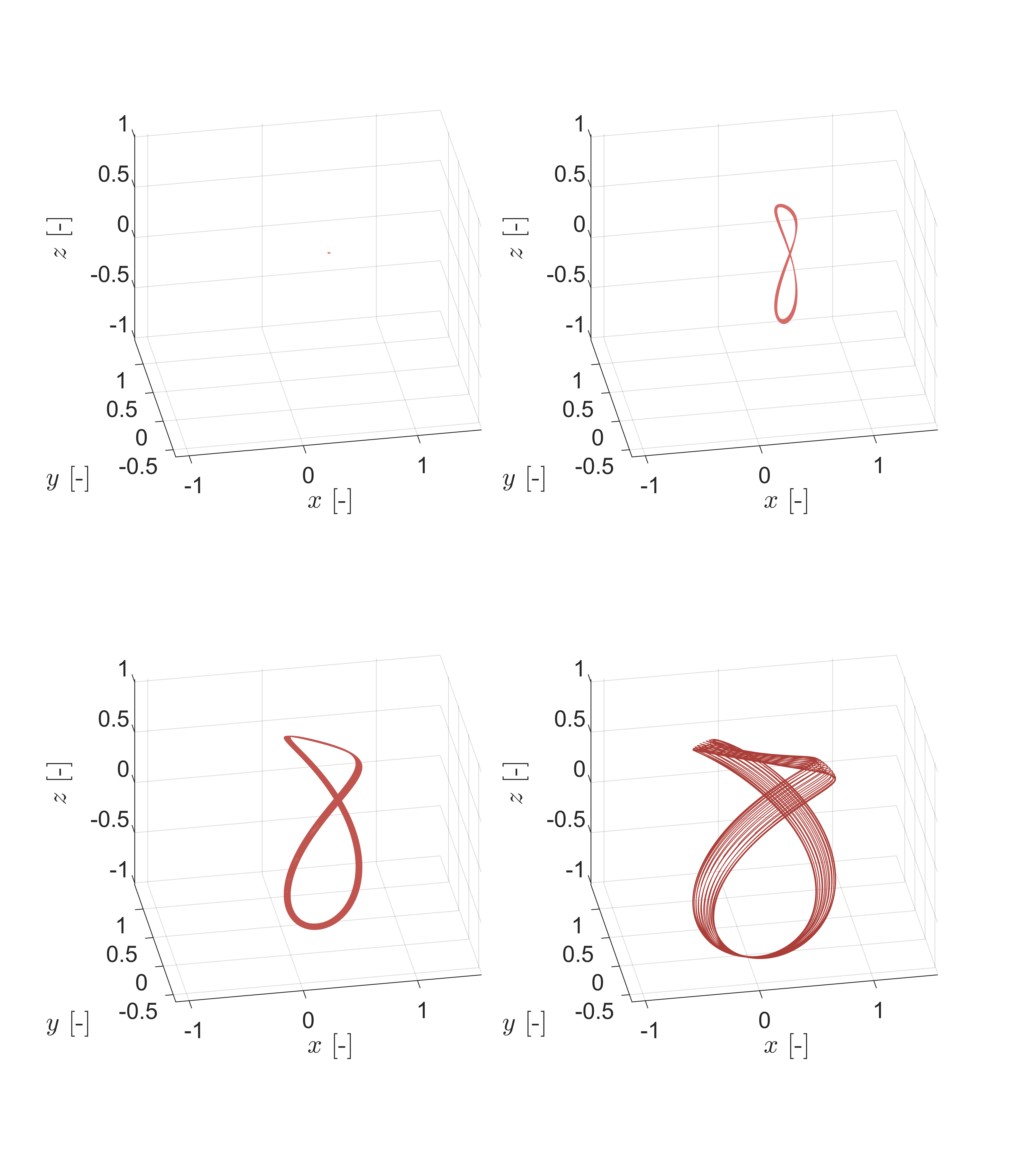}
     \end{subfigure}
        \caption{Hodograph and sample invariant tori of DE $L_4$ V family. See text for details.}
        \label{fig:DEL4_V}
\end{figure}

Figure \ref{fig:DEL4_V} shows a hodograph of the family computation, along with sample torus representations. Note that the family grows as $\omega_2$ increases. The invariant tori of this family retain a figure-8 shape as in the periodic vertical Lyapunov orbits about $L_4$ in the CR3BP; however, due to the periodic forcing of the Sun, an additional frequency is added. 

\begin{figure}
    \centering
    \includegraphics[width=\linewidth]{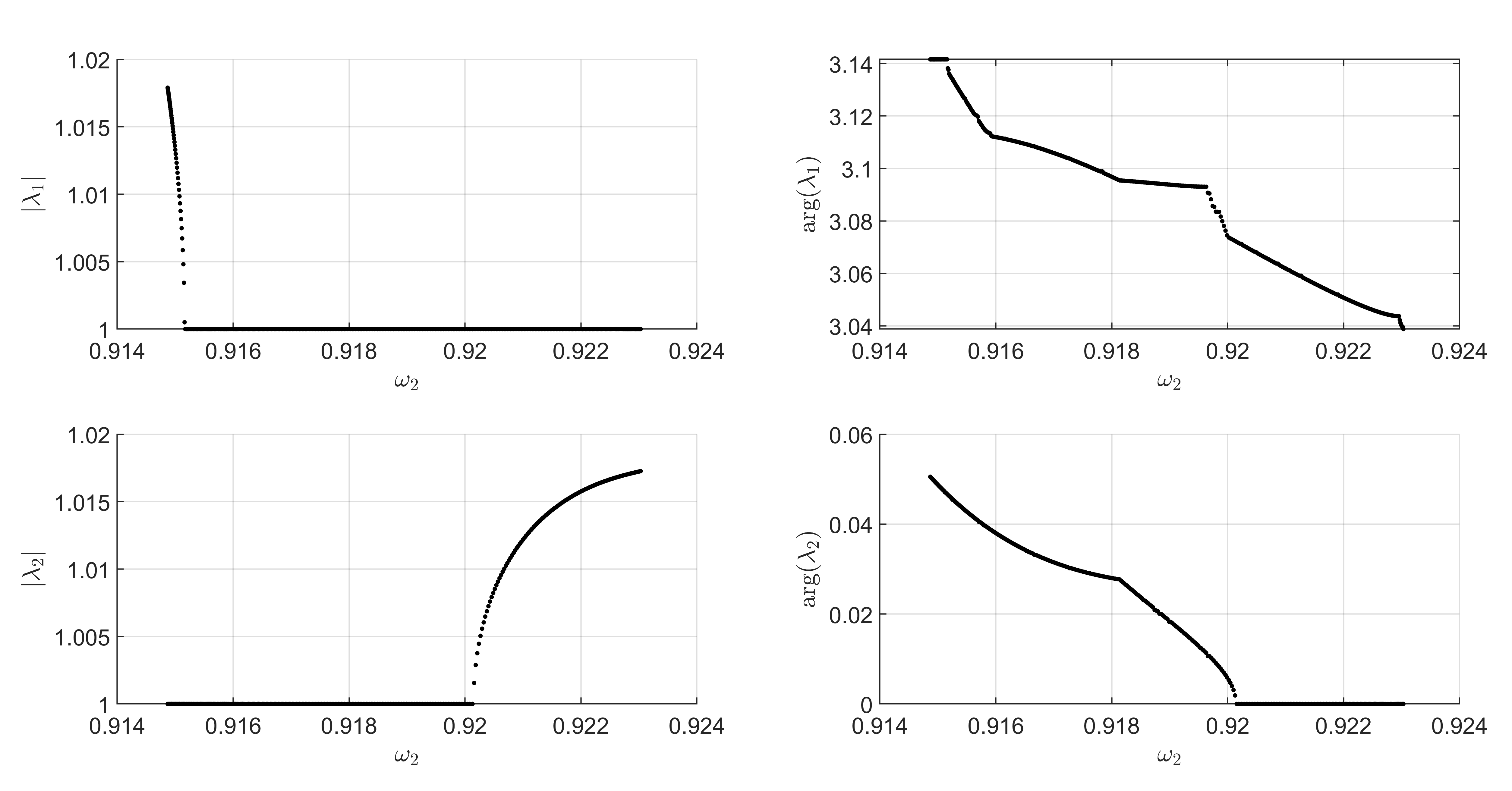}
    \caption{Normal behavior of DE $L_4$ V family of invariant tori.}
    \label{fig:Stab_DEL4_V}
\end{figure}

Figure \ref{fig:Stab_DEL4_V} shows the computed normal behavior of the family of tori. Unlike the DE $L_4$ H family, we show the normal behavior of DE $L_4$ V as a function of $\omega_2$ because the family is one-to-one in this variable. Once again, note that the family grows as $\omega_2$ grows. This figure shows that there are two bifurcations detected along the orbit family: a frequency-halving bifurcation, and a tangent bifurcation. The tangent bifurcation is seen in the CR3BP within the vertical Lyapunov family of periodic orbits around $L_4$, and so we say that it is an artifact of the CR3BP and leave its computation for future work, if desired. The frequency-halving bifurcation is not an artifact of the CR3BP and is thus born from the periodic forcing of the Sun. The term ``frequency-halving'' is a generalization of ``period-doubling'' from periodic orbits to quasi-periodic orbits (which have no period but instead have multiple frequencies). We note that the frequency-halving bifurcation is detected because a Krein collision of eigenvalues occurs along the negative real axis. 

We can summarize the normal behavior as follows. The family begins with the stability type of DE $L_4$: saddle $\times$ center. There is a frequency-halving bifurcation that changes the stability of DE $L_4$ V to type center $\times$ center. Then, there is a tangent bifurcation which changes the stability to center $\times$ saddle. 


\subsubsection[L4 2:1 H1]{$L_4$ 2:1 H1}
As the 2:1 resonant periodic orbit around $L_4$ is elliptic, there are two planar oscillatory modes. We adopt the naming convention of H1 and H2 for these horizontal families of invariant 2-tori. The $L_4$ 2:1 H1 family is born from the planar oscillation corresponding to the ``long-period'' (i.e., smaller frequency) central mode. Thus, the H2 family is born similarly but from the ``short-period'' (i.e., larger frequency) central mode. 

\begin{figure}[!htbp]
     \centering
     \begin{subfigure}[b]{0.95\textwidth}
         \centering
         \includegraphics[width=\textwidth]{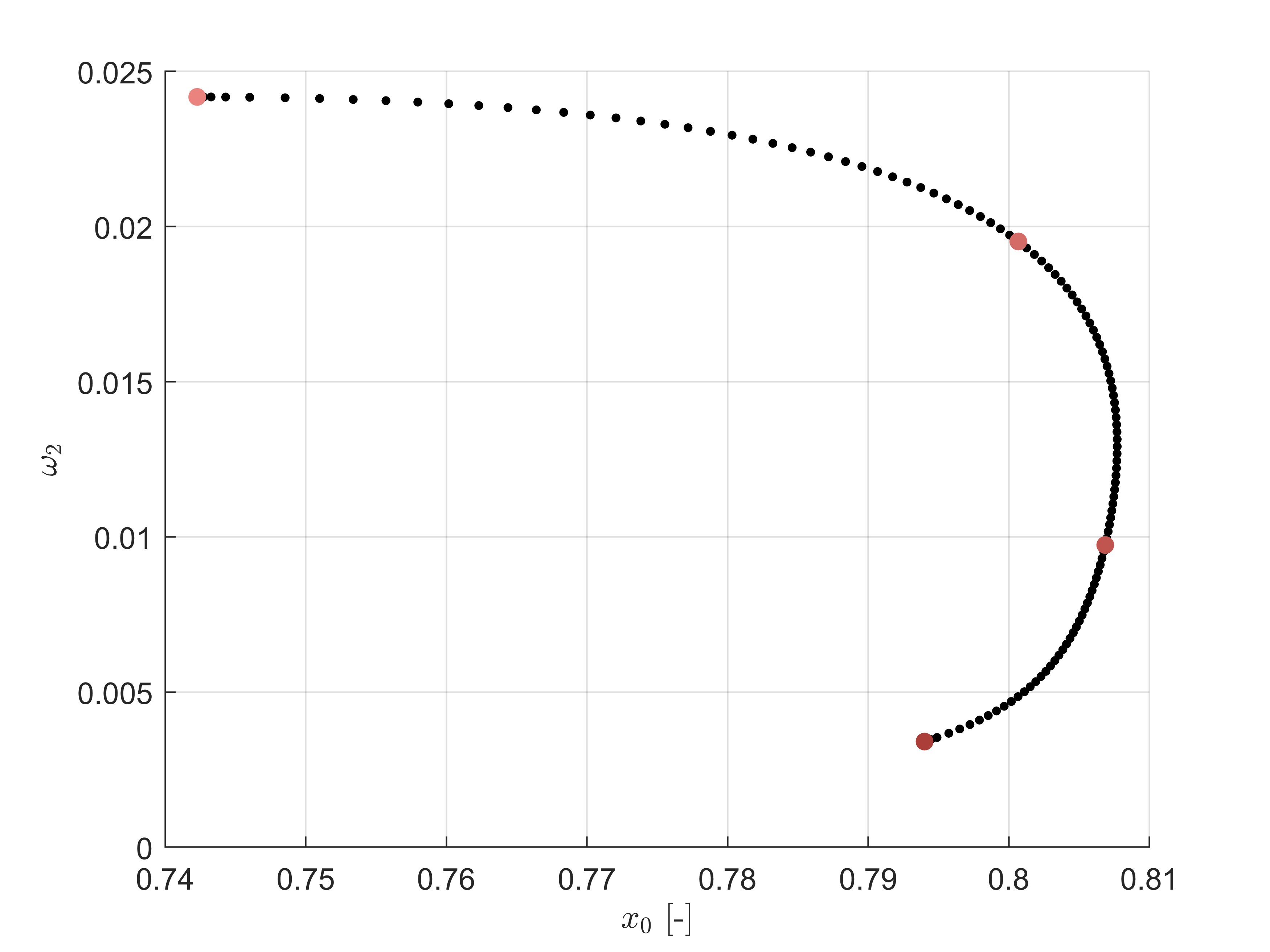}
     \end{subfigure}
     \hfill
     \begin{subfigure}[b]{0.95\textwidth}
         \centering
         \includegraphics[width=\textwidth]{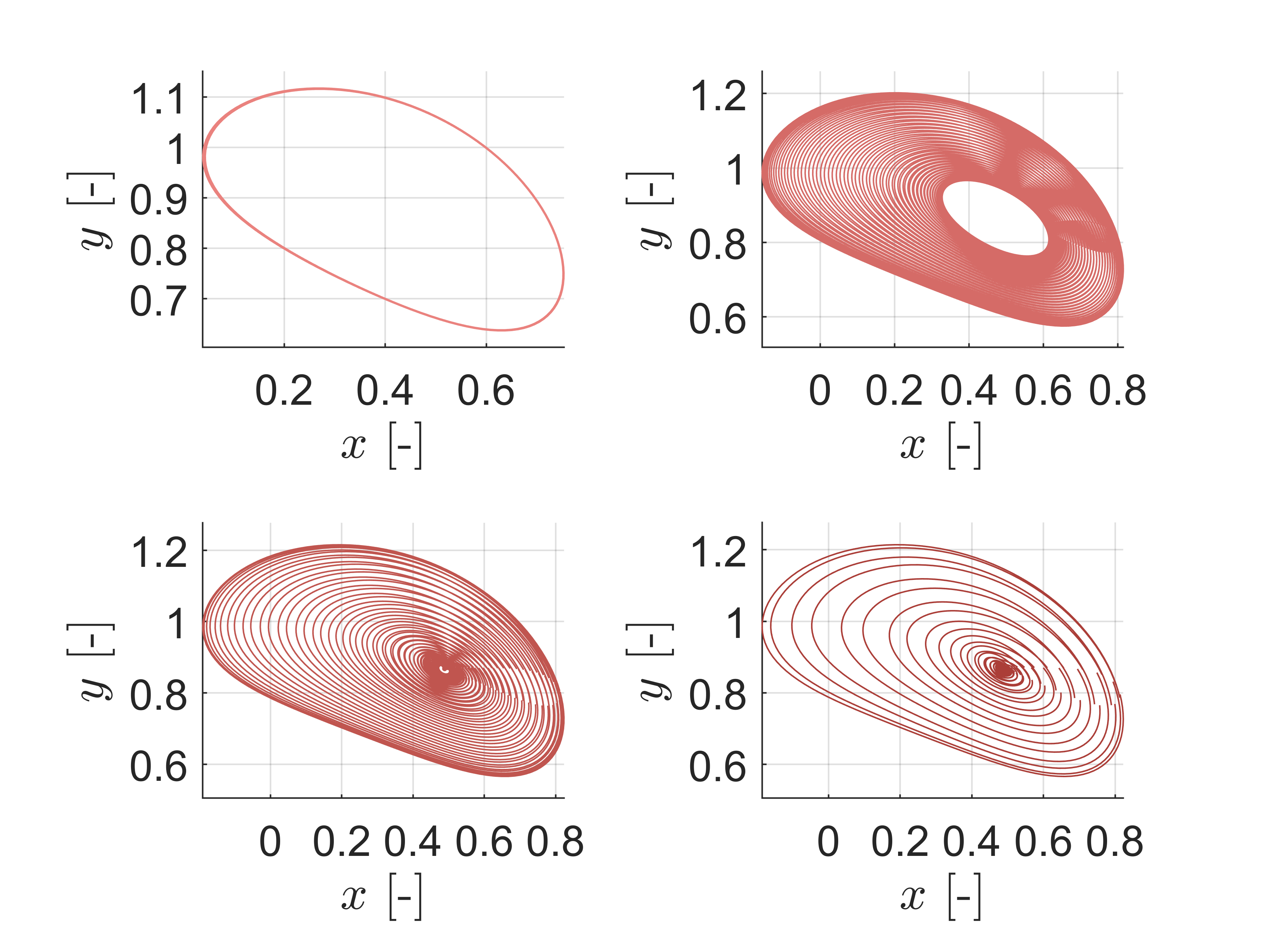}
     \end{subfigure}
        \caption{Hodograph and sample invariant tori of $L_4$ 2:1 H1 family. See text for details.}
        \label{fig:2to1_H1}
\end{figure}

\begin{figure}
    \centering
    \includegraphics[width=0.9\linewidth]{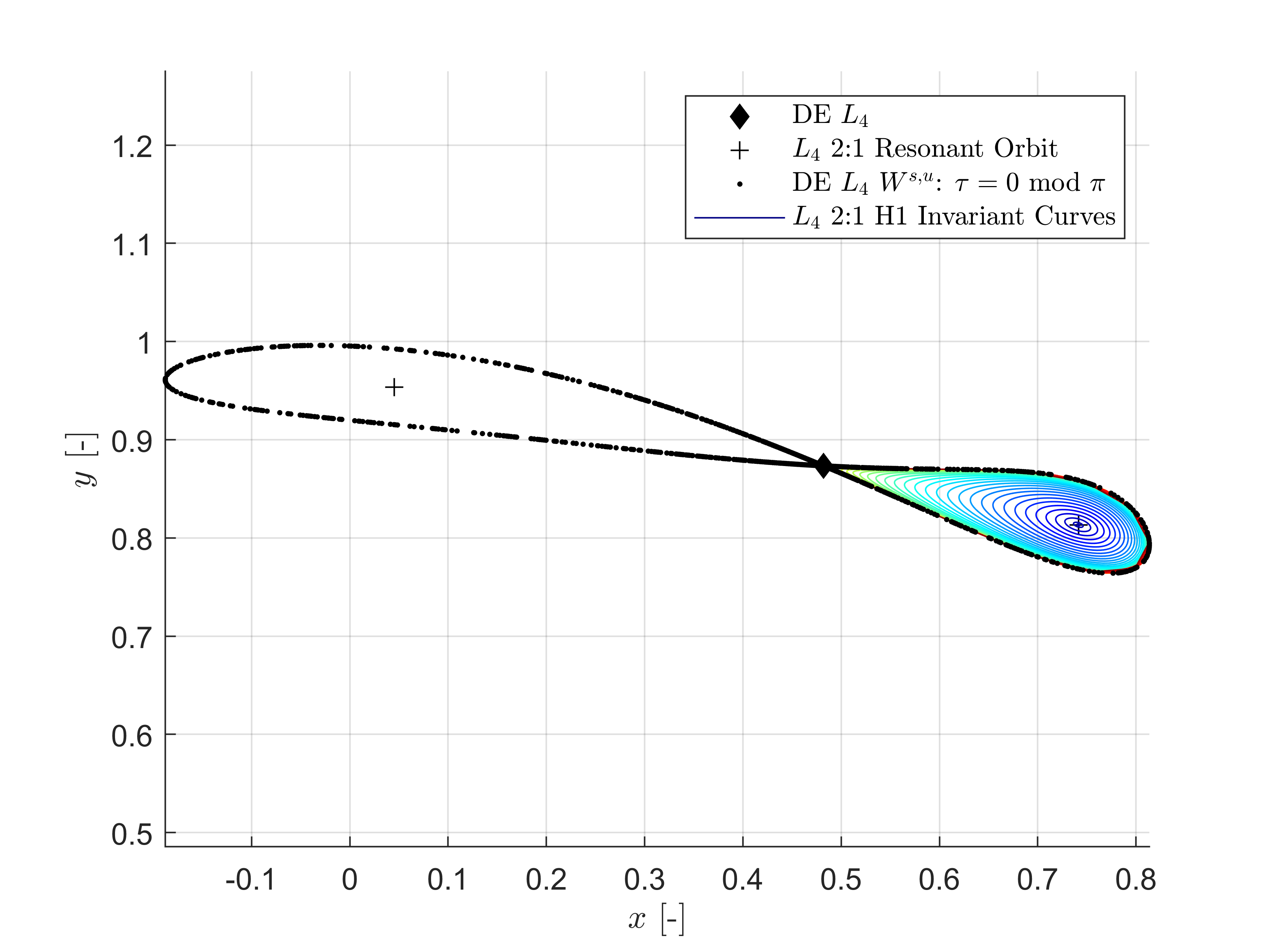}
    \caption{Invariant curves of the $\pi$-stroboscopic map with $W^{s,u}$ of DE $L_4$.}
    \label{fig:2to1_H1_ICs_Wsu}
\end{figure}

Figure \ref{fig:2to1_H1} shows a hodograph of the family computation, along with sample torus representations. Note that the family grows as $\omega_2$ decreases. Notice that the orbit appears to be getting pulled into DE $L_4$. While it has been shown that the stable and unstable manifolds of DE $L_4$ do not intersect \cite{scheeres1998restricted}, i.e., DE $L_4$ is not a homoclinic point, one can apply the following co-dimension argument to see that we expect there to be homoclinic and heteroclinic 2-tori near DE $L_4$. Consider an arbitrarily small planar 2-torus in the DE $L_4$ H family. Then the stable and unstable manifolds of this 2-torus intersect at a surface of dimension equal to the sum of the dimensions of the center-stable and center-unstable manifolds minus the dimension of the phase space. The center-stable and center-unstable manifolds each have dimension 3. Since we are considering the planar problem, the phase space is 5-dimensional (4 spatial, 1 temporal). Thus, $3 + 3 - 5 = 1$, so we expect there to be a 1-parameter family of homoclinic connections in the planar problem. These homoclinic connections between invariant 2-tori in the plane act as partial barriers to transport. It is clear from the stroboscopic map in Figure \ref{fig:2to1_H1_ICs_Wsu} that the invariant curves are approaching a boundary close to the stable and unstable manifolds of DE $L_4$. The authors speculate that the 2:1 H1 family ends in a homoclinic bifurcation with a homoclinic connection of an arbitrarily small member of the DE $L_4$ H family. The explicit computation of the connection is left for future work, as an identical argument shows 1-parameter families of heteroclinic connections exist, which may enable Arnold diffusion around DE $L_4$ in the Sun-Earth-Moon HR4BP. 

\begin{figure}[!htbp]
    \centering
    \includegraphics[width=\linewidth]{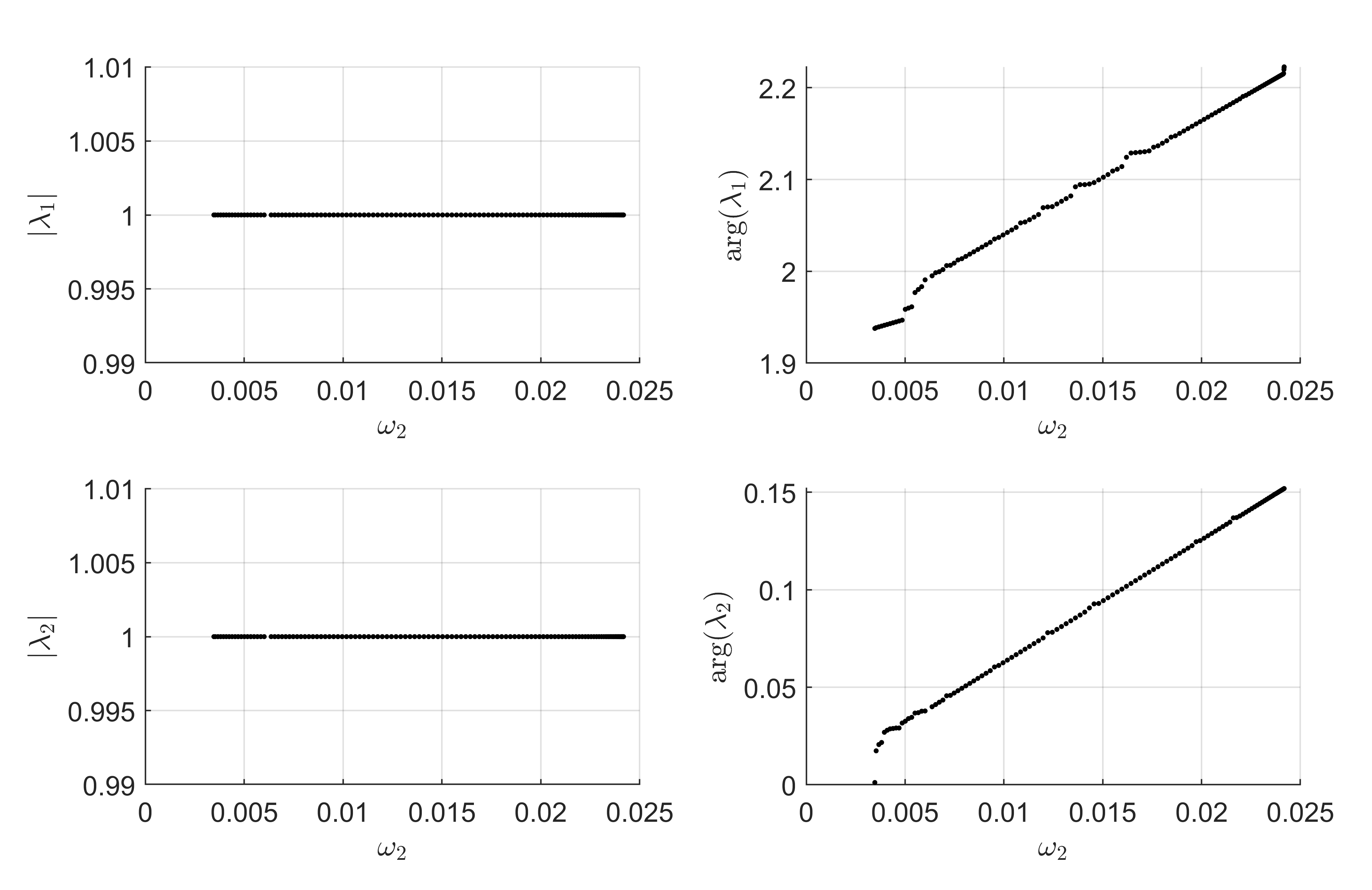}
    \caption{Normal behavior of $L_4$ 2:1 H1 family of invariant tori.}
    \label{fig:Stab_L4_2to1_H1}
\end{figure}

Figure \ref{fig:Stab_L4_2to1_H1} shows the computed normal behavior of the family of tori. As in the previous figure, the family grows in decreasing $\omega_2$. Observe that this family is elliptic, i.e., has normal stability type center $\times$ center. Yet, as we continue along this family of elliptic tori, by analyzing the bottom-right plot showing arg$(\lambda_2)$, we see that one pair of eigenvalues is approaching a Krein collision on the positive real axis, signaling an upcoming bifurcation. While we were unable to continue the family any further to exactly detect the bifurcation, this observation is consistent with the arguments made in the previous paragraph.

\subsubsection[L4 2:1 H2]{$L_4$ 2:1 H2}
\begin{figure}
     \centering
     \begin{subfigure}[b]{0.9\textwidth}
         \centering
         \includegraphics[width=\textwidth]{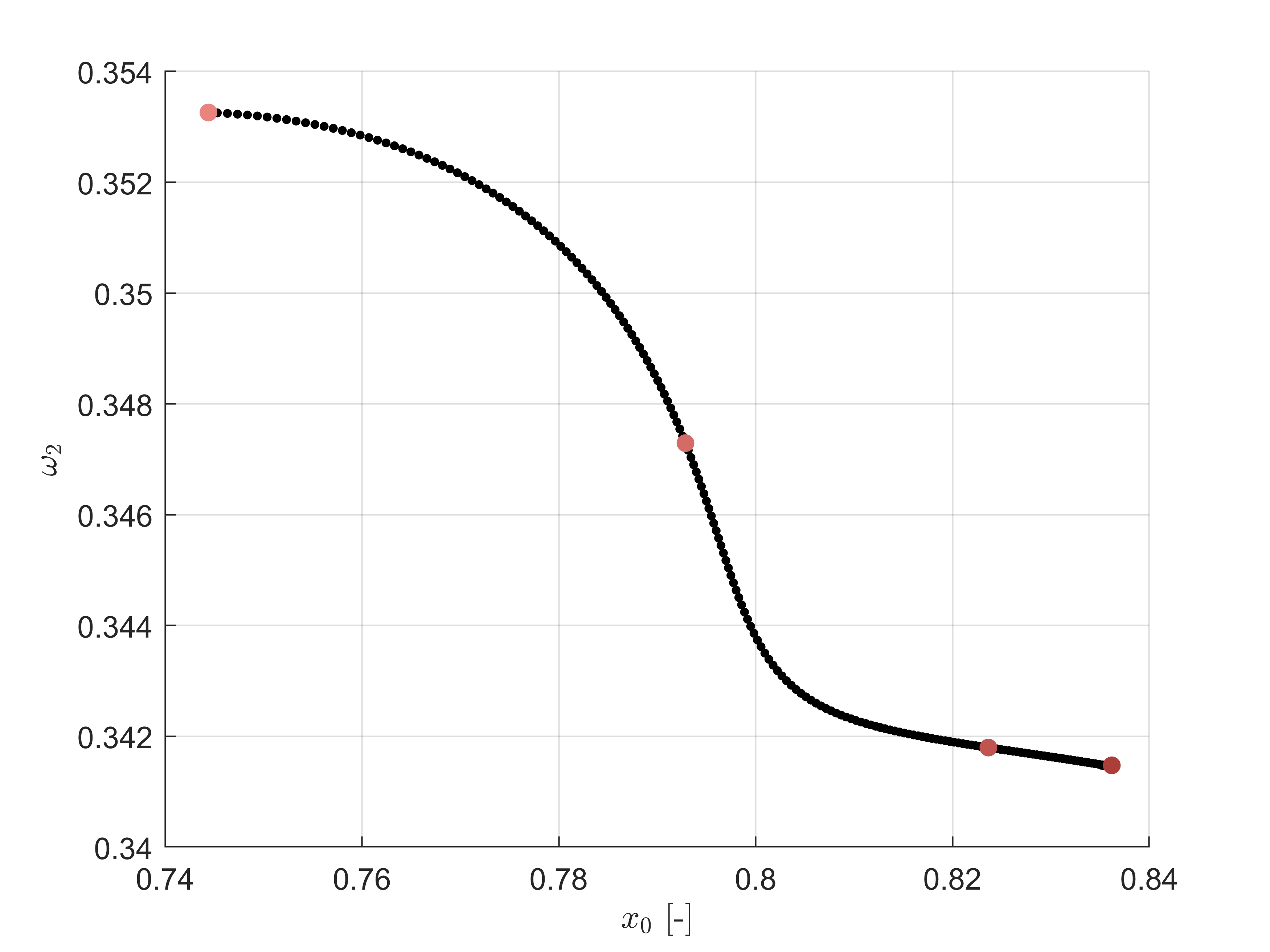}
     \end{subfigure}
     \hfill
     \begin{subfigure}[b]{0.9\textwidth}
         \centering
         \includegraphics[width=\textwidth]{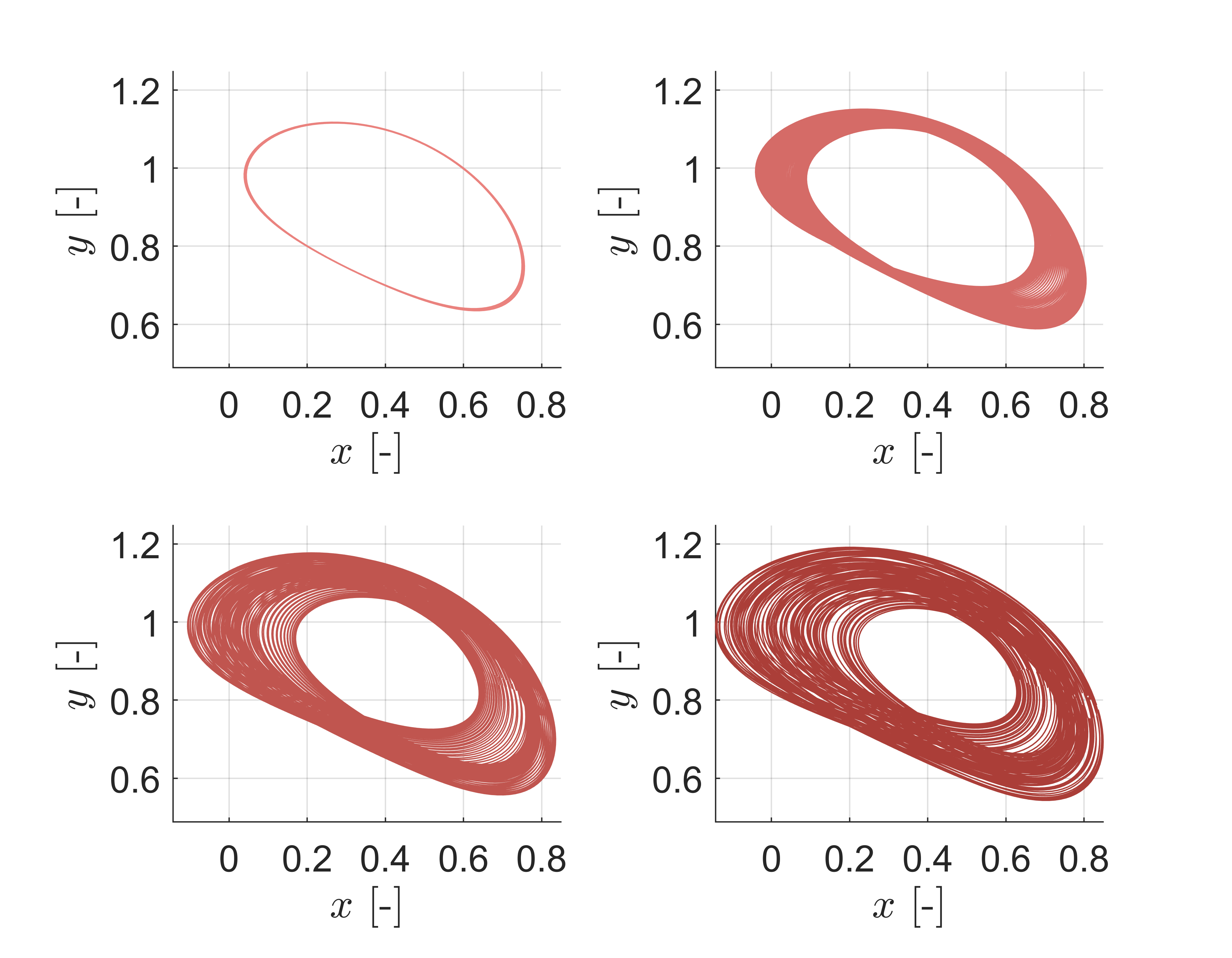}
     \end{subfigure}
        \caption{Hodograph and sample invariant tori of $L_4$ 2:1 H2 family. See text for details.}
        \label{fig:2to1_H2}
\end{figure}

\begin{figure}
    \centering
    \includegraphics[width=\linewidth]{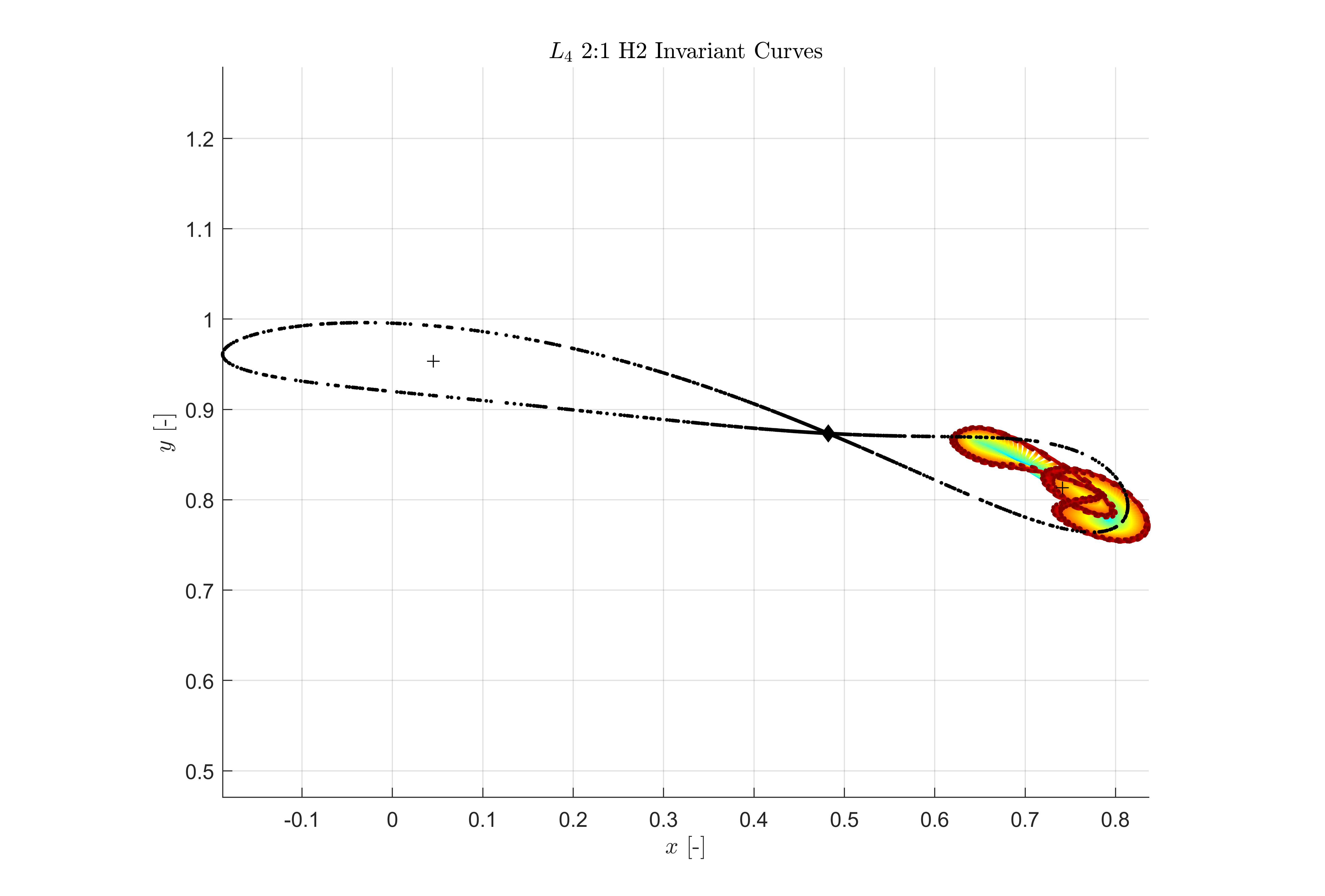}
    \caption{Invariant curves of the $\pi$-stroboscopic map. Note that $W^{s,u}$ of DE $L_4$ are not absolute barriers of transport.}
    \label{fig:2to1_H2_InvariantCurvesWsu}
\end{figure}

As stated in the previous section, the $L_4$ 2:1 H2 family of Lyapunov invariant 2-tori corresponds to the ``short-period'' (i.e., larger frequency) central mode of the elliptic 2:1 resonant periodic orbit. Figure \ref{fig:2to1_H2} shows a hodograph of the family computation, along with sample torus representations. Note that the family grows as $\omega_2$ decreases (or as $x_0$ increases). Unlike the 2:1 H1 family, the 2:1 H2 family is not pulled in by DE $L_4$. Instead, the tori in this family appear to similarly wind around the 2:1 periodic orbit. We can see that the homoclinic and heteroclinic connections of the DE $L_4$ H family do not impose such transport restrictions on the H2 family as in the previous section with the H1 family. Figure \ref{fig:2to1_H2_InvariantCurvesWsu} shows the invariant curves of H2 under the stroboscopic map. Note that as the colors move from blue to green to red, the family is increasing. Ultimately the computation of tori is stopped by the need to increase the number of Fourier modes above the maximum allowed in our computations--this threshold is reached due to the twisting of the invariant curves observed. 

\begin{figure}[!htbp]
    \centering
    \includegraphics[width=\linewidth]{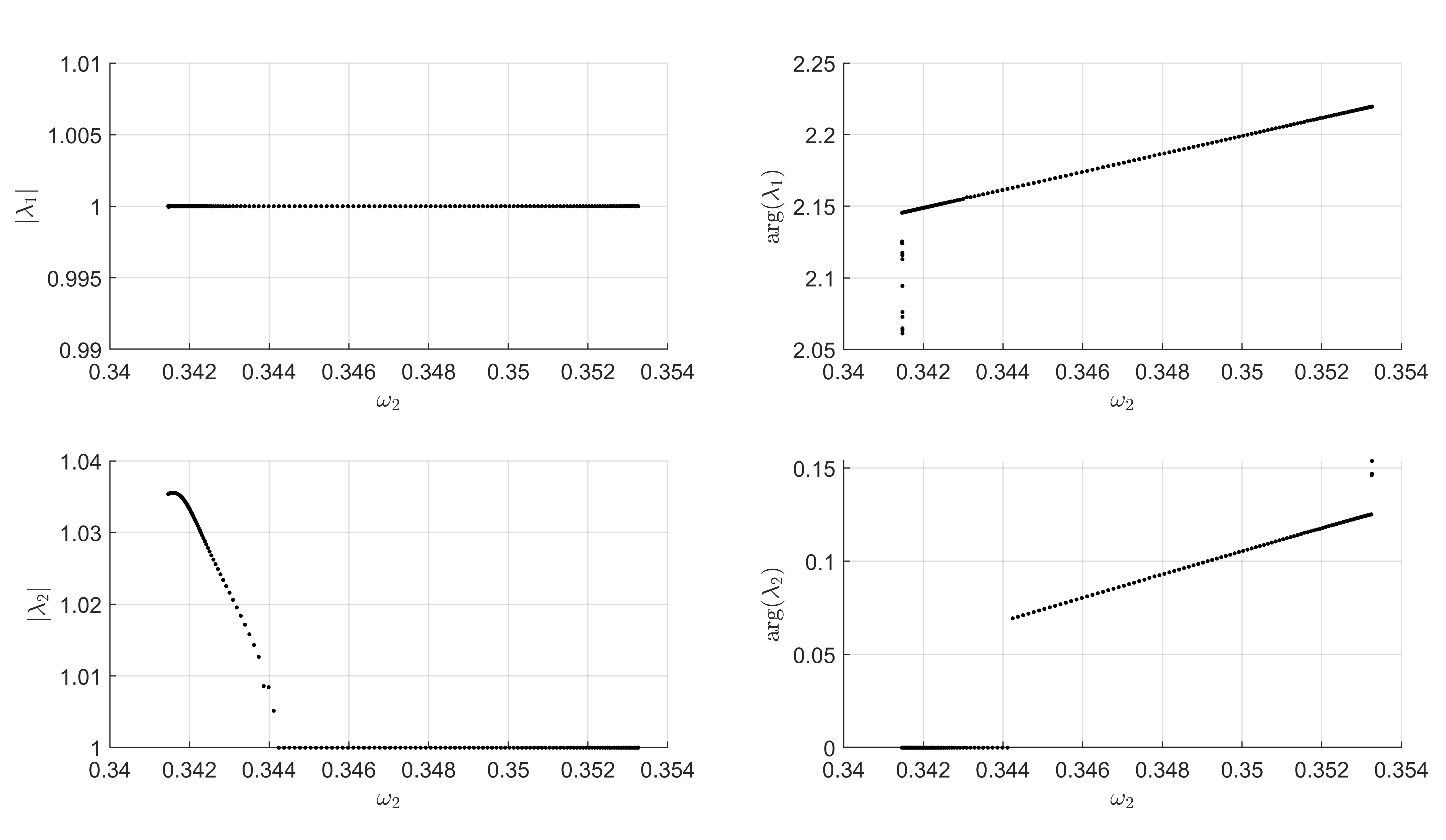}
    \caption{Normal behavior of $L_4$ 2:1 H2 family of invariant tori.}
    \label{fig:Stab_L4_2to1_H2}
\end{figure}

Figure \ref{fig:Stab_L4_2to1_H2} shows the computed normal behavior of the family of tori. As in the previous figure, the family grows in decreasing $\omega_2$. The family of invariant 2-tori inherits the normal behavior of the elliptic 2:1 resonant periodic orbit, i.e., the family begins with stability type center $\times$ center. There is one bifurcation identified in the family: a tangent bifurcation. This bifurcation connects the 2:1 H2 and V families, as can be seen in Figure \ref{fig:2to1_V} wherein one of the two broken branches bifurcates into the plane.

\subsubsection[L4 2:1 V]{$L_4$ 2:1 V}
\begin{figure}
    \centering
    \includegraphics[width=0.9\linewidth]{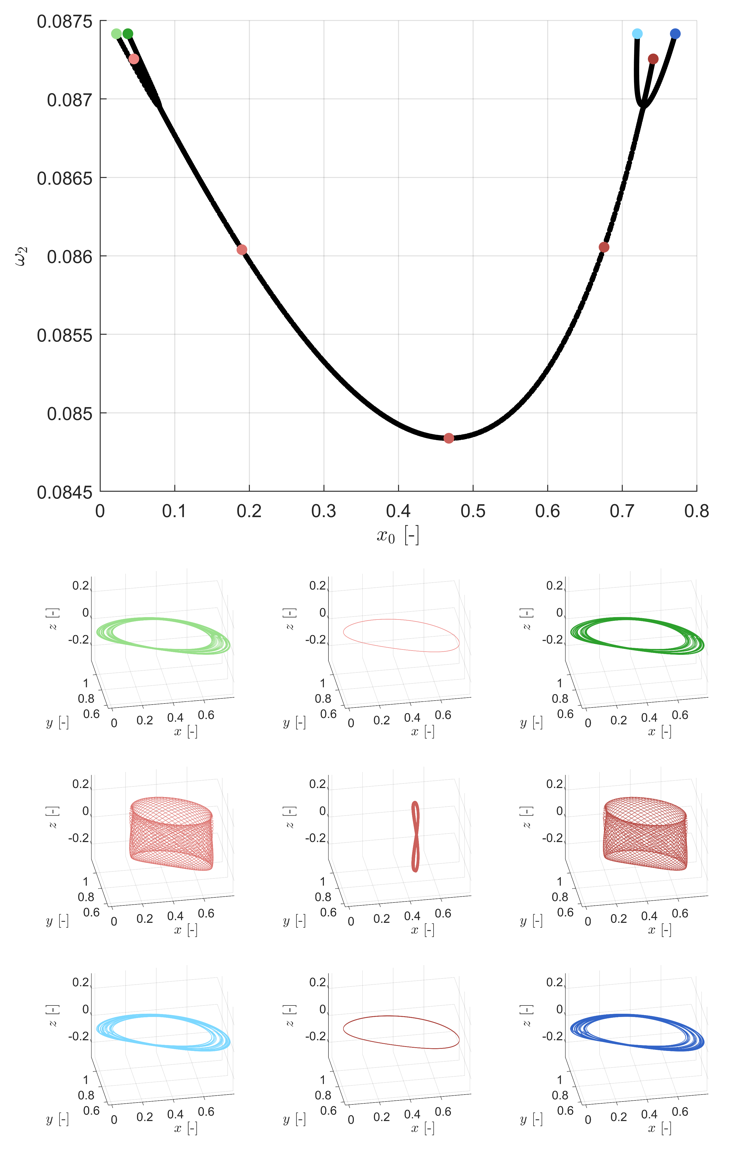}
    \caption{$L_4$ 2:1 V family with broken bifurcations. Note that there is a symmetry due to the 2:1 underlying period.}
    \label{fig:2to1_V}
\end{figure}

The final of the five families of the Lyapunov 2-tori computed near $L_4$ is the 2:1 V family. Figure \ref{fig:2to1_V} shows a hodograph of this family with sample torus representations. First, observe that there is a symmetry from the left and right sides of the hodograph; because these tori emanate from a 2:1 periodic orbit, showing the $\tau = 0 \mod \pi$ stroboscopic map yields two points representing the same torus on the hodograph. Moreover, the same principle applies to the previous two families, we simply omit one side because the two sides never meet, as in the 2:1 V family shown here. The family is computed from light to dark red. Due to the symmetry, the light blue and light green are equivalent representations of the same invariant torus (similarly from dark blue and dark green). The computation is stopped once the 2:1 periodic orbit is reached again. 

Figure \ref{fig:2to1_V} shows a bifurcation in the 2:1 V family, and the blue and green torus representations indicate that these bifurcations connect the planar 2:1 H2 family, as seen in the previous section. It is important to note that this bifurcation is broken and that the light and dark blues thus represent different invariant tori (similarly with the light and dark greens, by symmetry). The precise cause of this bifurcation is unknown to the authors. The normal behavior was computed to be elliptic for all 2-tori computed, including the tori born from the bifurcation. The stability is not presented because, as the tori were all found to be elliptic, the equivalence class representatives were not chosen--no clear additional insight would be gained by their inclusion. 

The orbit in the 2:1 V family which minimizes $\omega_2$ (the torus represented in the middle subplot of Figure \ref{fig:2to1_V}) is qualitatively similar to the DE $L_4$ V tori nearby--they are both in the shape of a figure 8. The frequency-halving bifurcation that occurs in the DE $L_4$ V family potentially relates the DE $L_4$ V and 2:1 V families of invariant tori. 

\subsection[Transport near L4 in the HR4BP]{Transport near $L_4$ in the HR4BP}\label{sec:Transport}
\begin{figure}[H]
    \centering
    \includegraphics[width=0.8\linewidth]{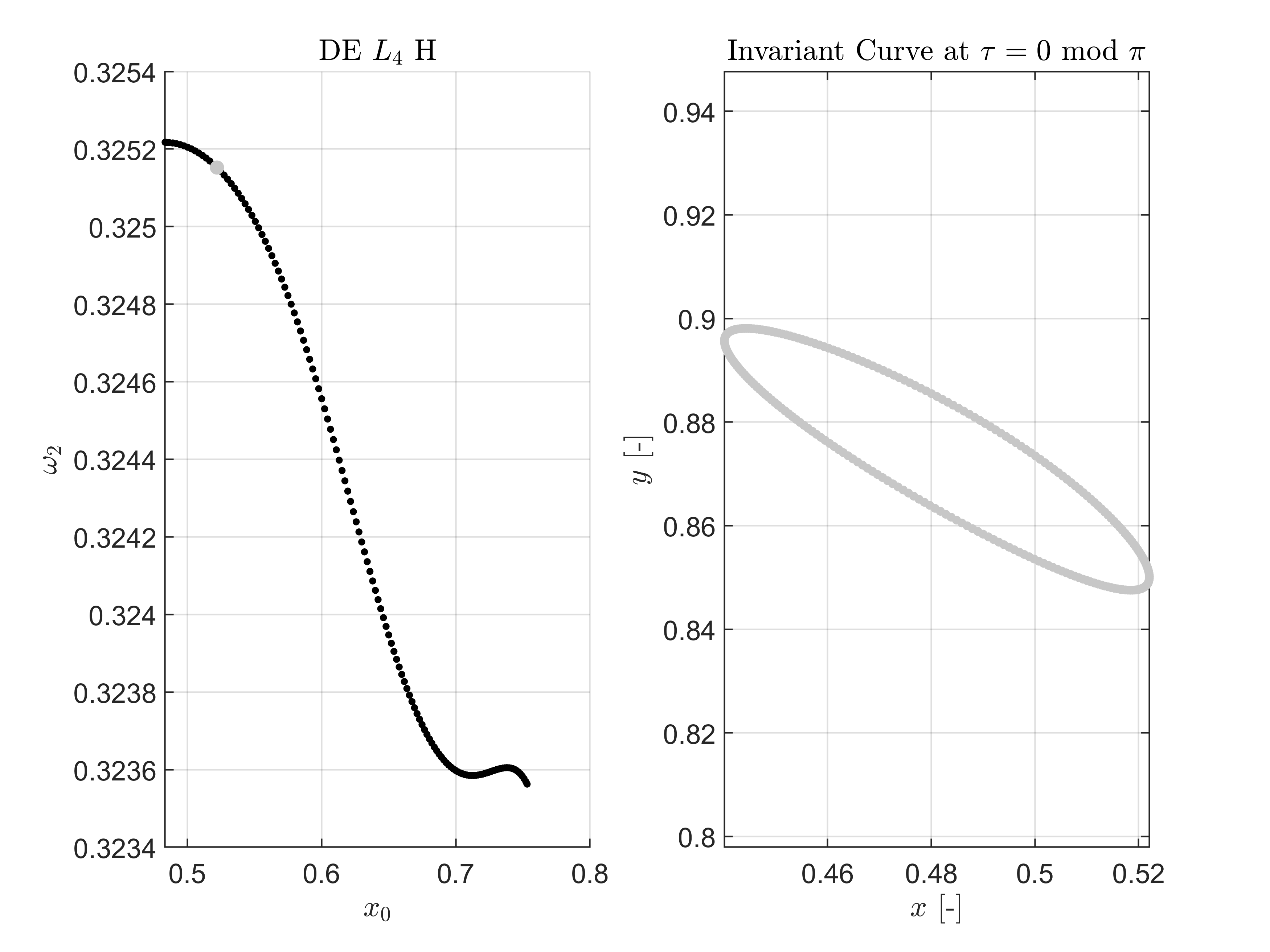}
    \caption{Selected invariant curve from DE $L_4$ H member.}
    \label{fig:DEL4H_IC}
\end{figure}

\begin{figure}[H]
    \centering
    \includegraphics[width=\linewidth]{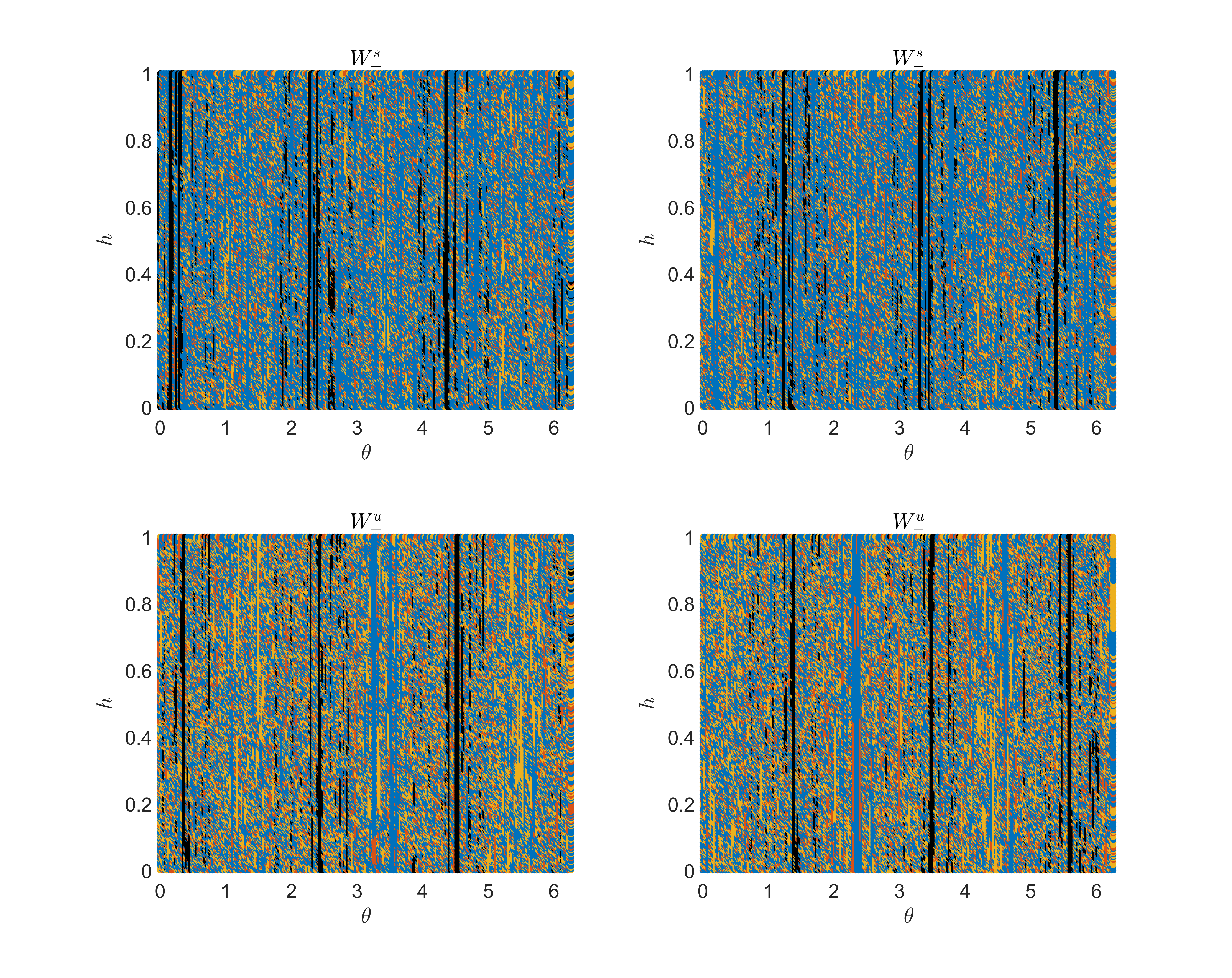}
    \caption{Sample stable manifold from a quasi-periodic invariant torus of the DE $L_4$ H family crossing GEO.}
    \label{fig:Transport_Maps}
\end{figure}

Two of the central questions about the dynamics around EM $L_4$ in the Sun-Earth-Moon system are: can a point start close to DE $L_4$ and escape the system? How are trajectory fates dispersed? In this section, we show that the hyperbolic invariant manifolds of an invariant torus near EM $L_4$ can flow close to Earth and Moon, as well as escape the system. As this is qualitatively different behavior from the Earth-Moon CR3BP, we choose an quasi-periodic orbit near DE $L_4$. In particular, we select a member of DE $L_4$ H, as shown in Figure \ref{fig:DEL4H_IC}. We take as the initial invariant curve to be at $\tau = 0$ mod $\pi$, and we use $h_0 = 10^{-5}$ for the initial perturbation onto the stable or unstable manifolds. Figure \ref{fig:Transport_Maps} shows the results of cylinder set computations around an invariant torus of DE $L_4$ H close to DE $L_4$. Blue points correspond to trajectories that escape the system; red points intersect the radius of GEO in the configuration space; yellow points impact the lunar surface; black points do not satisfy any of the above criteria. The statistics for the cylinder sets are given in the table below:

\begin{center}
    \begin{table}[!htbp]
    \begin{tabular}{ccccc}
    \hline
    Manifold & Lunar Impact (\%) & GEO Radius (\%) & Escapes (\%) & None (\%) \\
    \hline
    $W_+^s$ & 17.33 & 13.24 & 60.75 & 8.69 \\
    $W_-^s$ & 16.89 & 13.69 & 61.42 & 8.00 \\
    $W_+^u$ & 27.78 & 16.52 & 48.25 & 7.45 \\
    $W_-^u$ & 27.29 & 17.23 & 48.64 & 6.85 \\
    \hline
    \end{tabular}
    \end{table}
\end{center}

\begin{figure}[H]
    \centering
    \includegraphics[width=0.7\linewidth]{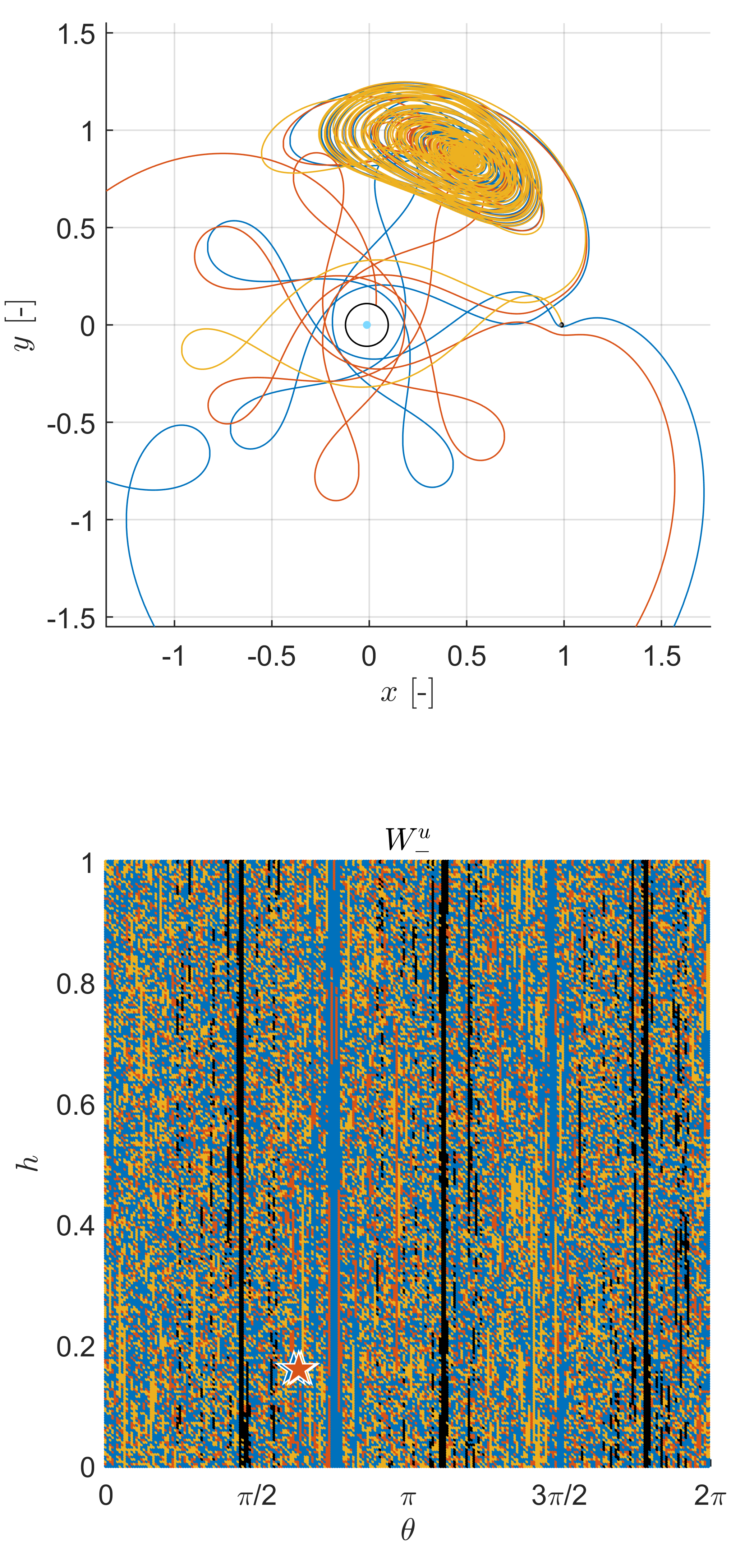}
    \caption{Trajectory comparison between adjacent points of $W_{-}^u$ map with different fates.}
    \label{fig:TransportComparison}
\end{figure}

The chaotic distribution of fates on each cylinder set can be described via closest approach to the Moon. We take as an illustrative example in the $W_{-}^u$ map, choosing three points within a small neighborhood with each a distinct fate. Figure \ref{fig:TransportComparison} shows the resulting trajectories (top) with the corresponding points on the cylinder (bottom). Note that the trajectories remain close together until the first lunar perilune, whence the trajectories enter the interior region near the Earth. Now, the yellow trajectory intersects the lunar surface. The red trajectory spends more time in the interior region then has a second close approach to the Moon and temporarily orbits the system before intersecting the GEO radius. Finally, the blue trajectory similarly winds around Earth in the interior region before making a close approach with the Moon--notably closer than the red trajectory--and picks up significant energy from the flyby causing it to escape the system. 

\begin{figure}[H]
    \centering
    \includegraphics[width=\linewidth]{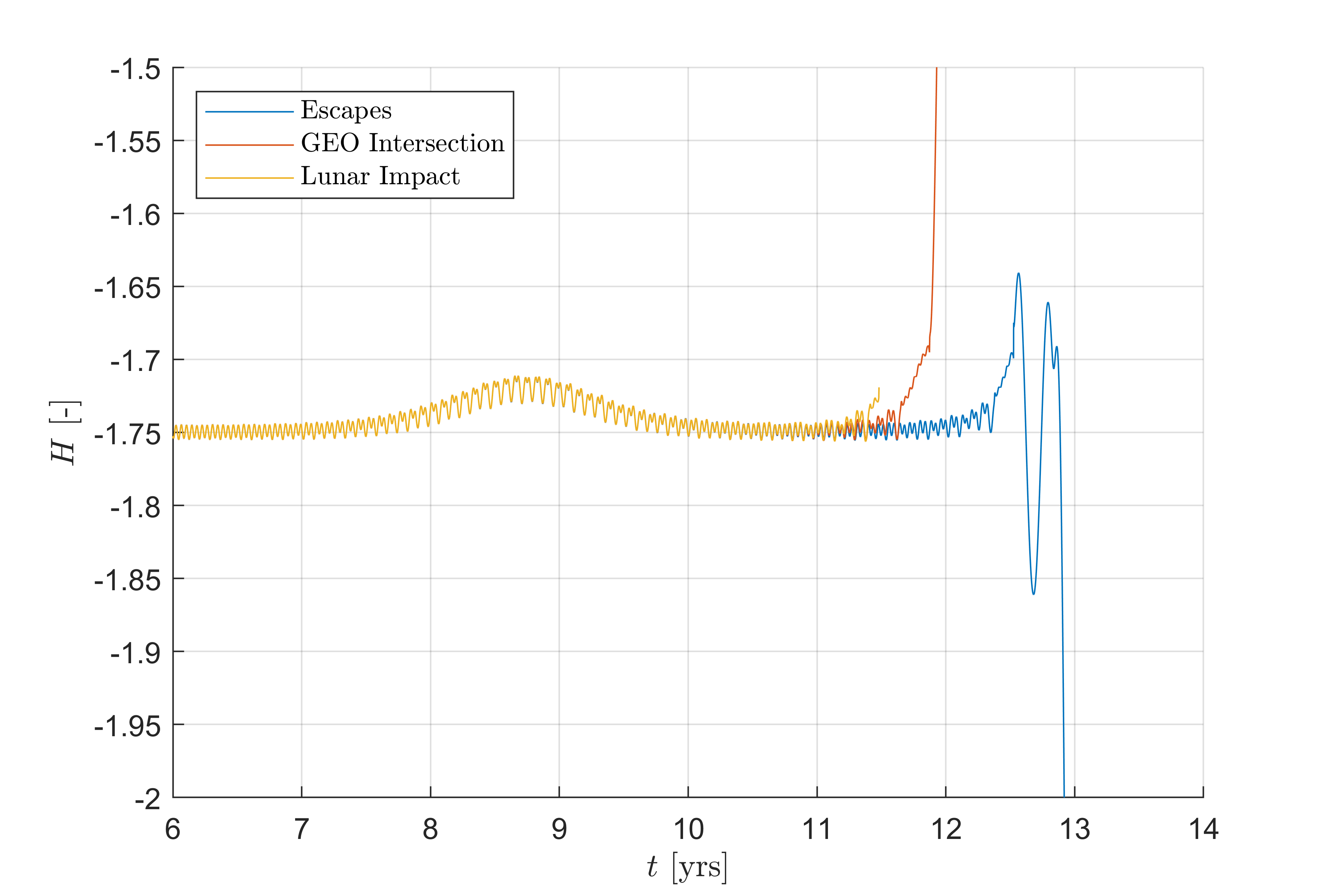}
    \caption{Energy comparison between adjacent map points with different fates.}
    \label{fig:TransportComparison_Energy}
\end{figure}

\begin{figure}[H]
    \centering
    \includegraphics[width=\linewidth]{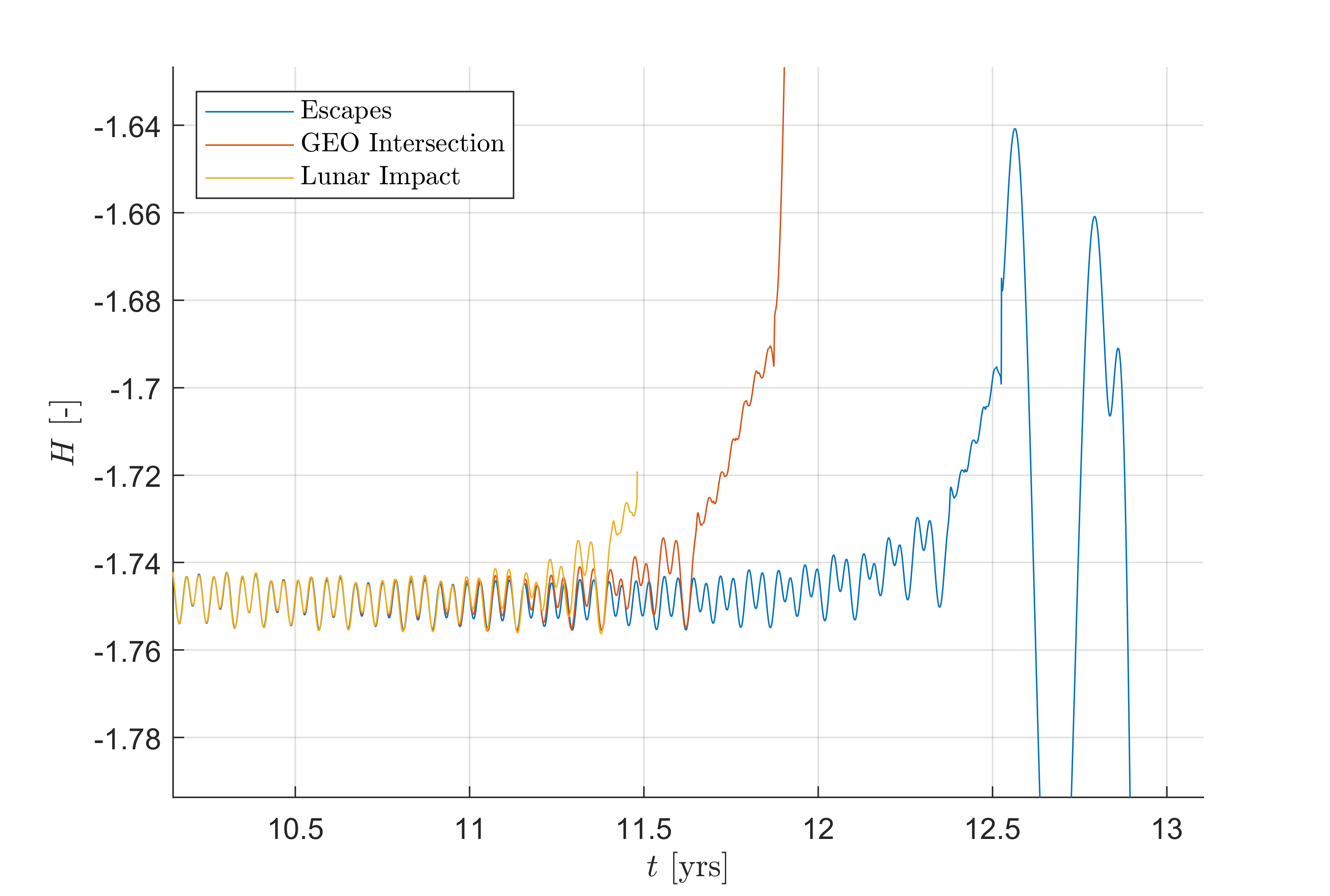}
    \caption{Zoom of energy comparison between adjacent map points with different fates.}
    \label{fig:TransportComparison_EnergyZoom}
\end{figure}

The change in energy from each lunar flyby can be seen by computing the Hamiltonian along each trajectory, as shown in Figures \ref{fig:TransportComparison_Energy} and \ref{fig:TransportComparison_EnergyZoom}. The Hamiltonian is not constant along the flow as the system is periodic. Figure \ref{fig:TransportComparison_Energy} shows the difference in energy change across each flyby, especially red and blue, i.e., GEO intersection and escapes the system. In the zoomed plot in Figure \ref{fig:TransportComparison_EnergyZoom}, we can more clearly see the effect of the first lunar perilune on the energy of each trajectory. Namely, the difference in energy along each trajectory is matched until the first lunar perilune, whence the energies are dephased and diverge from each other. Clearly, the yellow trajectory diverges least, as it intersects the lunar surface and hence has no lunar flyby. 

\begin{figure}[H]
    \centering
    \includegraphics[width=0.8\linewidth]{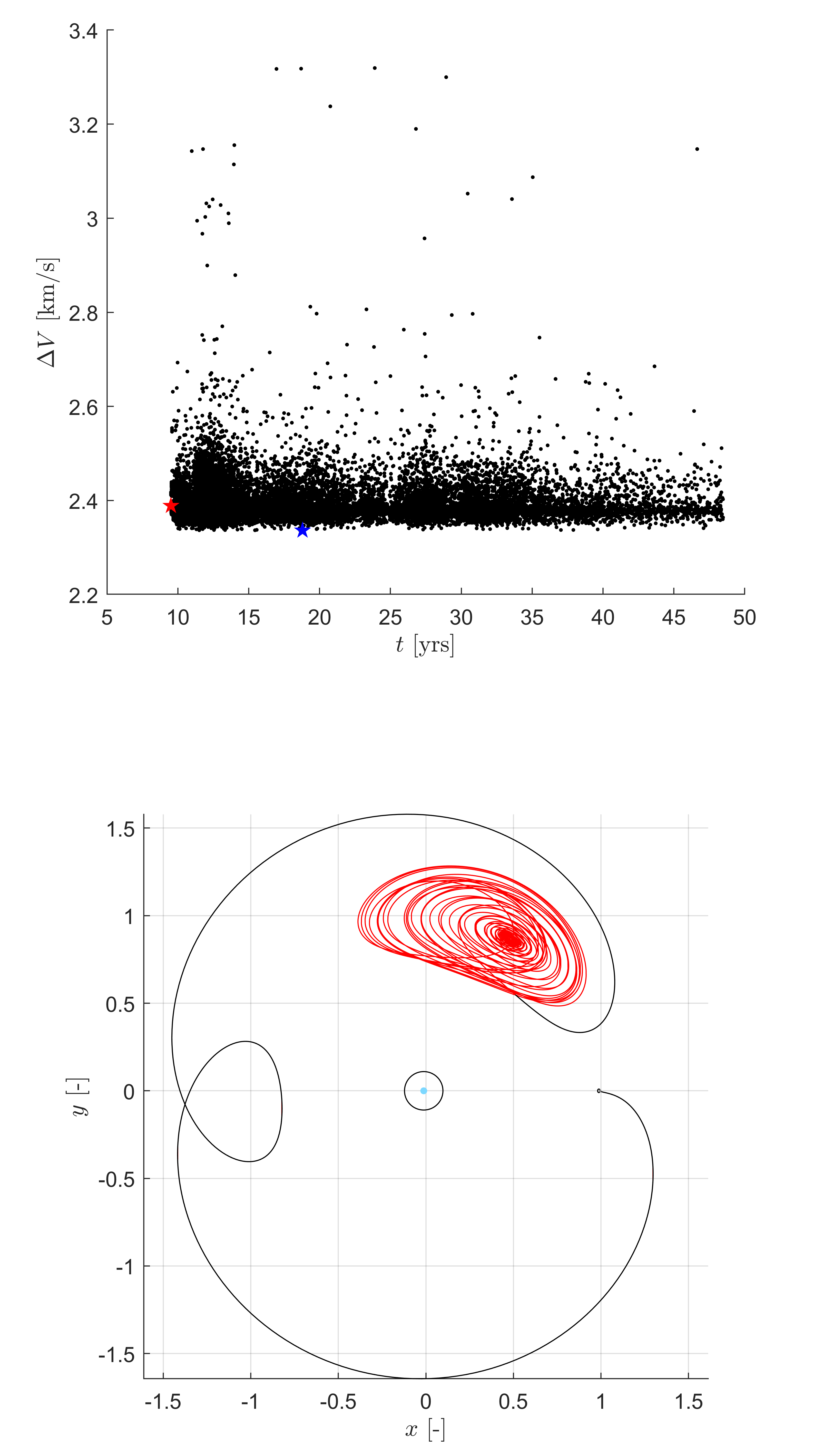}
    \caption{$\Delta V$ [km/s] and time of flight [yrs] for direct transfers from the lunar surface onto the stable manifold of the selected small-amplitude quasi-periodic orbit near $L_4$, including minimum time of flight trajectory. See text for details.}
    \label{fig:Transport_Pareto}
\end{figure}

Finally, as it is possible to move around the Earth-Moon system via stable and unstable manifolds of a quasi-periodic orbit close to EM $L_4$, we can consider $L_4$ as a region of interest for the disposal of objects near the Moon. As we have seen, the stable manifold of the selected quasi-periodic orbit intersects the lunar surface at many points. We can compute the necessary $\Delta V$ for a single-impulsive maneuver sending a particle from the lunar surface directly onto the stable manifold of the chosen orbit, accounting for the velocity gained from the rotation of the Moon. This gives a conservative estimate on disposal maneuver cost for spacecraft along periodic and quasi-periodic orbits near $L_1$ and $L_2$, as there is less energy difference between these orbits to apply a maneuver compared with the lunar surface. Figure \ref{fig:Transport_Pareto} shows the comparison of all points along the stable manifold, $W_{\pm}^s$, in time of flight (in years) and $\Delta V$ (in km/s), as well as the trajectory corresponding to the minimum time of flight. In the top plot, the red star is the minimum time of flight, and the blue star is the minimum $\Delta V$. In the bottom plot, the first segment, shown in black, is the portion of the trajectory transiting into the greater $L_4$ region, while the second segment, shown in red, remains in the greater $L_4$ region. While the time of flight to the specific quasi-periodic orbit is nearly 10 years, the transit time to the greater $L_4$ region is approximately 78 days. 

\section{Comparison with Existing Models}\label{sec:Comparison}
To gain a more complete picture of the dynamics near EM $L_4$, we compare the findings presented in this work with existing work in similar restricted 4-body problem models. This is especially imperative to place our work in the context of the literature. In this section, we compare the periodic orbits near EM $L_4$ of the HR4BP with those of the bicircular restricted 4-body problem (BCP) and the quasi-bicircular restricted 4-body problem (QBCP). Most literature describing dynamics in the Sun-perturbed Earth-Moon system takes the BCP as the dynamical model. In particular, dynamics near EM $L_4$ have been studied to varying extents by Jorba et al. in the BCP, QBCP, and a JPL model. Hence, we also make remarks about numerical simulations done in the JPL ephemeris model based on $L_4$ orbits in the BCP.

\subsection{Bicircular Restricted 4-Body Problem}\label{sec:BCP}
The bicircular restricted 4-body problem (BCP) is a dynamical system that can be seen as a $2\pi$-periodic perturbation of the circular restricted 3-body problem (CR3BP) that accounts for the direct effect of the Sun. The BCP is a Hamiltonian system described by Hamiltonian function, defining the corresponding momenta as $p_x = \dot{x} - y$, $p_y = \dot{y} + x$, and $p_z = \dot{z}$:
\begin{equation}
H_{\text{BCP}} = \frac{1}{2} \left( p_x^2 + p_y^2 + p_z^2 \right) + y p_x - x p_y - \frac{1-\mu}{r_{\text{PE}}} - \frac{\mu}{r_{\text{PM}}} - \frac{m_S}{r_{\text{PS}}} - \frac{m_S}{a_S^2}\left( y \sin \theta - x \cos \theta \right),
\end{equation}
where $r_{\text{PE}}^2 = (x-\mu)^2 + y^2 + z^2$, $r_{\text{PM}}^2 = (x+1-\mu)^2 + y^2 + z^2$, $r_{\text{PS}}^2 = (x-x_S)^2 + (y-y_S)^2 + z^2$, $x_S = a_S \cos \theta$, $y_S = -a_S \sin \theta$, and $\omega_S = \omega_S t$. The frequency of the Sun, $\omega_S = 0.925195985518290$, plays an important role in the continuation of periodic solutions from the CR3BP into the BCP (compared to the HR4BP). Figure \ref{fig:BCP} shows a schematic of the BCP \cite{andreu1998quasi}. Note that in these two models, the positions of the Earth and Moon are rotated $180^\circ$ in the $x-y$ plane. 

Note that the indirect effect of the Sun, i.e., the effect of the Sun's gravity on the other two primary bodies, is not modeled in the BCP. By contrast, the HR4BP accounts for both the direct and indirect effects of the Sun. The quasi-bicircular problem (QBCP) is a modification of the BCP which incorporates the indirect effect of the Sun. Furthermore, a continuation parameter $\varepsilon \in \mathbb{R}$ is introduced to transition between CR3BP and BCP models as follows:
\begin{equation}
    H_{\text{BCP}}^\varepsilon = H_{\text{CR3BP}} + \varepsilon \widehat{H}_{\text{BCP}}, \quad \widehat{H}_{\text{BCP}} = -m_S \left( \frac{1}{r_{\text{PS}}} + \frac{y\sin\theta - x\cos\theta}{a_S^2} \right).
\end{equation}
From the above equation, one can see that $H_{\text{BCP}}^{\varepsilon = 0} = H_{\text{CR3BP}}$ and $H_{\text{BCP}}^{\varepsilon = 1} = H_{\text{BCP}}$. A similar process is applied in the QBCP. 

\begin{figure}
    \centering
    \includegraphics[width=\linewidth]{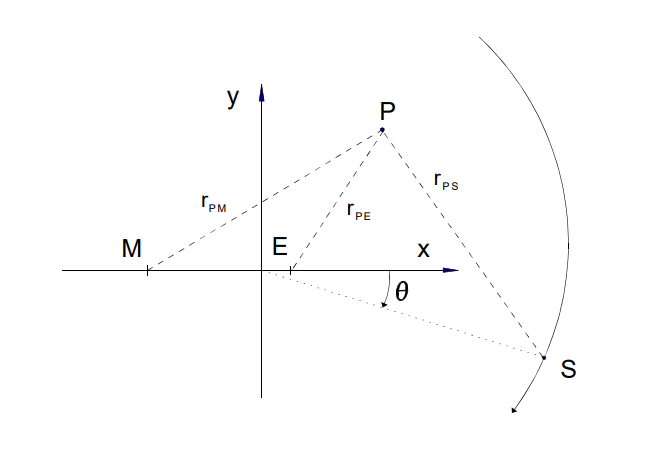}
    \caption{Schematic of Sun-Earth-Moon-particle motion in the bicircular problem (BCP) \cite{andreu1998quasi}.}
    \label{fig:BCP}
\end{figure}

\subsubsection{Periodic Solutions}

In the bicircular problem, as the linear frequencies around $L_4$ are different from $2k\pi\omega_S$ $\forall k \in \mathbb{Z}$, the $L_4$ periodic orbit can be continued to a $2\pi/\omega_S$-periodic orbit for $\varepsilon > 0$ small. Figure \ref{fig:BCP_L4} shows on the left plot the continuation diagram for $L_4$ in the BCP \cite{castella2000vertical}. Observe that there is a broken pitchfork bifurcation that causes a loss of uniqueness for the dynamical substitute of $L_4$ as $\varepsilon \to 1$. In the right plot of \ref{fig:BCP_L4}, the projection onto the configuration space is showed for the periodic orbits replacing $L_4$. The normal behavior for PO1 is center (horizontal) $\times$ saddle $\times$ center (vertical); both PO2 and PO3 are elliptic periodic orbits. Note that PO1 of the BCP and DE $L_4$ in the HR4BP share stability characteristics which is responsible for the instability in the region around $L_4$ in the Earth-Moon system. Similarly, note that PO2 and PO3 of the BCP share stability properties with the 2:1 resonant orbit of the HR4BP. Not only do they share stability aspects, but the amplitudes of PO2 and PO3 are nearly identical to the $L_4$ 2:1 orbit in the HR4BP--compare Figures \ref{fig:BCP_L4} and \ref{fig:L4_PeriodicOrbits}. In fact, the $L_4$ 2:1 orbit in the HR4BP shares the phasing properties of PO2 and PO3, which have opposite phasing; however, rather than two separate periodic orbits with opposite phasing, the $L_4$ 2:1 orbit has a period double that of the system allowing for the opposite phasing to occur on a single periodic orbit. 

Figure \ref{fig:BCP_QBCP_L4} shows a comparison between the continuation of $L_4$ in the BCP and in the QBCP. Due to their similarity in bifurcations and normal behavior, the literature considers the BCP as a sufficient approximation to model the dynamics around $L_4$ using the BCP or QBCP. 

One key difference between the HR4BP and QBCP is the period of the perturbation. On one hand, the BCP and QBCP are periodic with period $2\pi/\omega_S$, where $\omega_S \approx 0.9252$ is taken to be the frequency of the Sun. On the other hand, the HR4BP is periodic with period $\pi$ due to the symmetry provided by the Hill approximation of the Sun's gravitational effect. Consequently, the periodic orbit replacing $L_4$ in the BCP and QBCP will have period $2\pi/\omega_S$, rather than period $\pi$ as in the HR4BP. Similarly, the 2:1 resonant periodic orbit in the CR3BP is a different orbit when considering the periodic forcing of the QBCP or the HR4BP. In the HR4BP the 2:1 resonant orbit has period $2\pi$ in the CR3BP; in the QBCP the qualitatively similar resonant orbit would be 1:1 resonant orbit with period $2\pi/\omega_S$. This CR3BP periodic orbit was not able to be continued into the Sun-Earth-Moon BCP. Despite this, the behavior of periodic solutions around $L_4$ is qualitatively identical between the BCP and HR4BP.

\begin{figure}
    \centering
    \includegraphics[width=\linewidth]{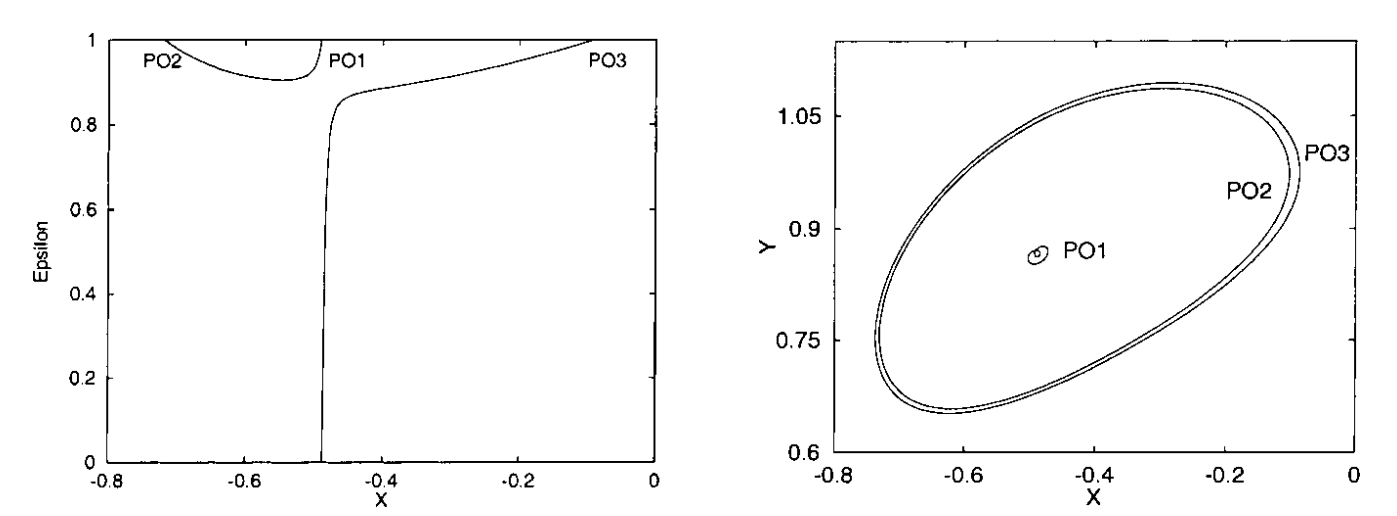}
    \caption{Continuation of $L_4$ in the bicircular problem (BCP) \cite{castella2000vertical}.}
    \label{fig:BCP_L4}
\end{figure}

\begin{figure}
    \centering
    \includegraphics[width=0.7\linewidth]{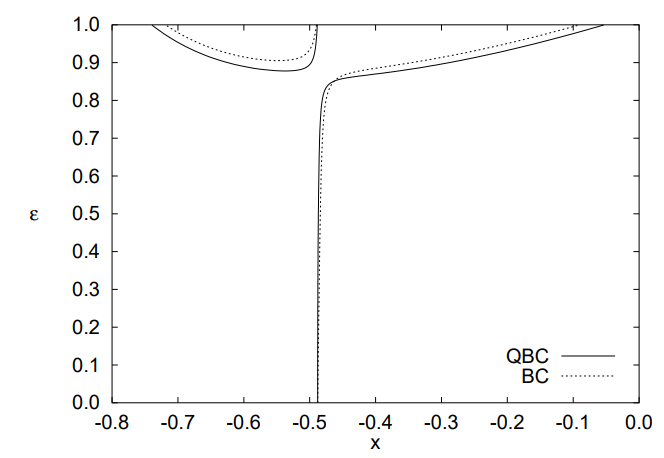}
    \caption{Continuation of $L_4$ in the quasi-bicircular problem (QBCP) compared to the BCP \cite{andreu1998quasi}.}
    \label{fig:BCP_QBCP_L4}
\end{figure}

\subsubsection{Invariant Tori}
In \cite{castella2000vertical}, the vertical families of invariant 2-tori around PO1, PO2, and PO3 are computed in the Sun-Earth-Moon BCP. The planar families of invariant 2-tori around $L_4$ in the BCP have not been presented in the literature, though \cite{gimeno2024effect} computes some planar families of Lyapunov 2-tori in the BCP with the perturbation of solar radiation pressure. 

Figure \ref{fig:BCP_L4V} shows the continuation of $L_4$ vertical families of invariant 2-tori in the BCP for several values of $\varepsilon$ \cite{castella2000vertical}. The horizontal axis shows the value of $\dot{z}$ of the invariant curve when $z = 0$, and the vertical axis shows the rotation number. In the Sun-Earth-Moon BCP, i.e., when $\varepsilon = 1$, there are 3 vertical families of invariant 2-tori emanating off of PO1, PO2, and PO3--these are denoted by F1, F2, and F3. Observe that there is a broken pitchfork bifurcation between the three families of invariant tori, similar to the bifurcation of periodic orbits. In fact, one can see from the continuation of these families in $\varepsilon$ that the bifurcation of periodic orbits is the genesis of the relationship between the families of tori. The relationship between vertical tori, especially the relationship between PO2/PO3 and PO1, is mirrored in the HR4BP; moreover, while the bifurcation in DE $L_4$ V has not been computed, preliminary calculations at smaller values of $m$ suggest there may be a similar link relating the DE $L_4$ V and $L_4$ 2:1 V families. The difference in the HR4BP is that the underlying bifurcation is not so clearly generated by the bifurcation of periodic orbits as in the BCP. This is under current investigation. 

\begin{figure}
    \centering
    \includegraphics[width=\linewidth]{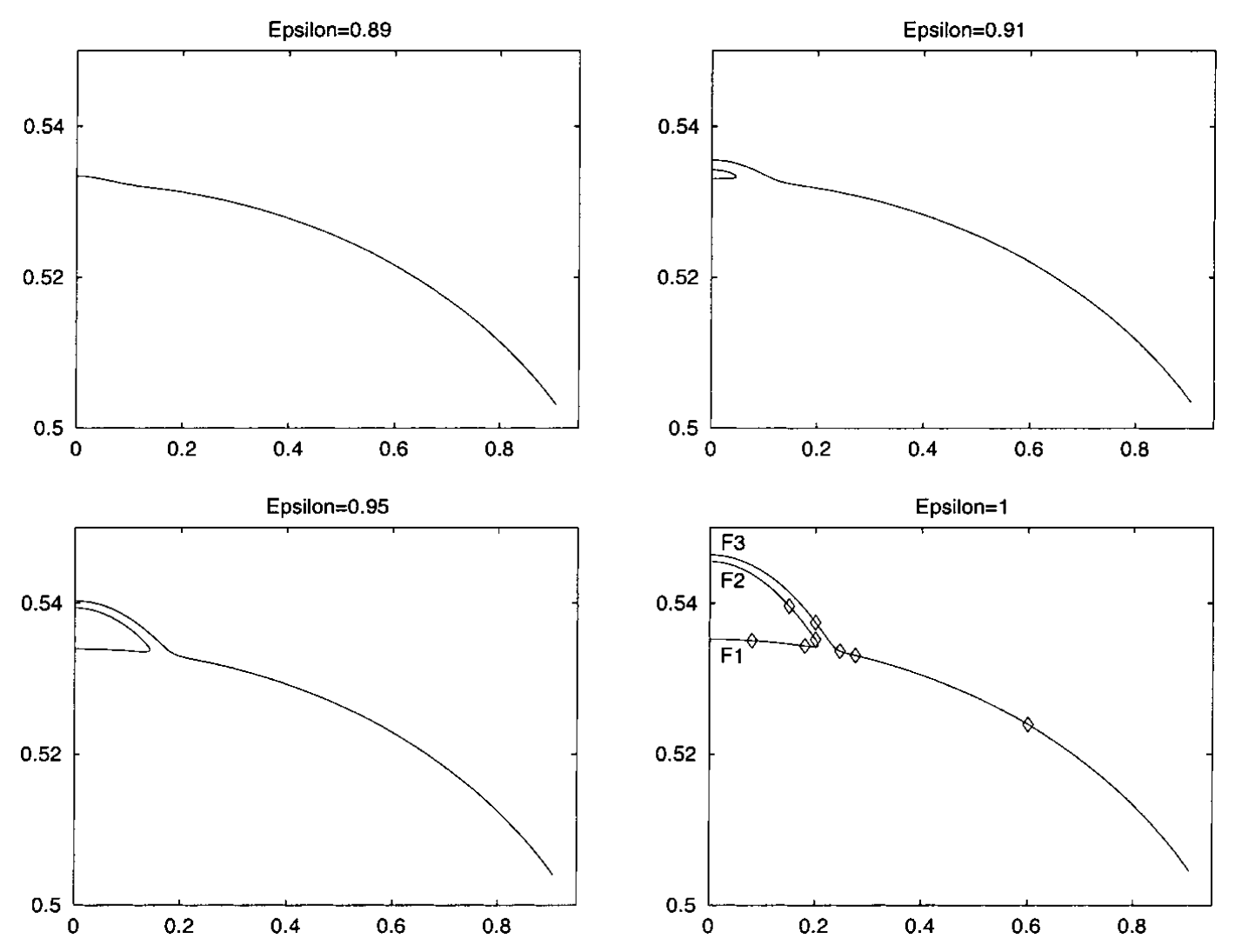}
    \caption{Continuation of $L_4$ vertical families of invariant 2-tori in the BCP for several values of $\epsilon$ \cite{castella2000vertical}. The horizontal axis shows the value of $\dot{z}$ of the invariant curve when $z = 0$, and the vertical axis shows the rotation number.}
    \label{fig:BCP_L4V}
\end{figure}

\subsubsection[Natural and Powered Connections to L4]{Powered Transfers to $L_4$}
Finally, as we have discussed transport phenomena around $L_4$ and the utility for spacecraft disposal in cislunar space, we put this work into the context of the limited related literature. While there are no works considering transfers near the Moon to $L_4$, there is some work considering transfers from Earth to $L_4$ in the BCP. Hence, a topic of continued study is to compute transfers between Earth and $L_4$, which could be used for a communications satellite relaying information between Earth and the lunar gateway around EM $L_2$.

Tan et al. consider single impulsive transfers from Earth using the window of easily approach (WOEA) without exploiting invariant manifolds--neither hyperbolic invariant manifolds nor quasi-periodic invariant tori \cite{tan2020single}. While the dynamics in the Sun-perturbed model allow for the exploitation of natural pathways along stable manifolds of invariant tori around DE $L_4$, Tan et al. fail to take advantage of these natural structures which can alleviate some of the fuel cost in $\Delta V$.

Conversely, Liang et al. exploit quasi-periodic orbits near $L_3$ as parking orbits to compute powered transfers from Earth to $L_4$ \cite{liang2021leveraging}. In their work, stable manifolds of several $L_3$ planar invariant tori are computed via backward propagation to find intersections with parking orbits around the Earth, ranging from low parking orbits to geostationary orbits. Unstable manifolds of the same $L_3$ planar invariant tori are then computed to identify intersections with stable regions around $L_4$ at a given epoch. A codimension argument shows that there are expected to be heteroclinic connections between 2-tori around $L_3$ and $L_4$, but this is not considered in their work.

\section{Conclusions}\label{sec:Conclusions}
In this paper we have investigated the motion of a small particle moving near the Earth-Moon triangular points using as our model the Hill restricted 4-body problem (HR4BP). The periodic forcing of the Sun changes the $L_4$ equilibrium point into an isolated $\pi$-periodic periodic orbit, called DE $L_4$. This perturbation qualitatively changes the linear behavior around the triangular points from elliptic to partially hyperbolic. Additionally, we show the continuation of the 2:1 resonant periodic orbit, explicitly computing the bifurcation killing one of the two branches that continue into the HR4BP. We considered the semi-analytical method of center manifold reduction, which we found to be extremely limited in utility. Hence, we computed 5 families of invariant 2-tori using the flow map method around DE $L_4$ and the 2:1 resonant periodic orbit, as well as their stability. A bifurcation relating to he vertical families of invariant 2-tori was identified. Finally, we investigated transport phenomena near EM $L_4$ by computing cylinder sets of a particular small planar quasi-periodic orbit close to DE $L_4$. We compared results to existing literature of similar work in the bicircular restricted 4-body problem. 

\backmatter

\bmhead{Acknowledgments}
The authors thank Damennick Henry for thoughtful discussions regarding the dynamical model and computation of invariant tori.



\bibliography{references}

\end{document}